\documentclass{article}
\usepackage{amssymb}
\usepackage{amsfonts}
\usepackage{amsfonts}
\usepackage{mathrsfs}
\usepackage {amsmath}
\usepackage {CJK}
\usepackage {amsthm}
\usepackage{graphicx}

\title{UNIQUENESS OF THE FOLIATION OF CONSTANT MEAN CURVATURE SPHERES IN ASYMPTOTICALLY FLAT 3-MANIFOLDS}
\author{Shiguang Ma}

\date{}

\theoremstyle{definition}
\newtheorem{lemma}{Lemma}[section]
\newtheorem{definition}[lemma]{Definition}
\newtheorem{theorem}[lemma]{Theorem}
\newtheorem{proposition}[lemma]{Proposition}
\newtheorem{corollary}[lemma]{Corollary}
\newtheorem{example}[lemma]{Example}
\newtheorem{remark}[lemma]{Remark}
\renewcommand{\proof}{Proof. }
\numberwithin{equation}{section}
\begin{document}
 \maketitle
\begin{abstract}
In this paper I study the constant mean curvature surface in
asymptotically flat 3-manifolds with general asymptotics. Under some
weak condition, I prove that outside some compact set in the
asymptotically flat 3-manifold with positive mass, the foliation of
stable spheres of constant mean curvature is unique.
\end{abstract}

\section{Introduction}
A three-manifold $M$ with a Riemannian metric $g$ and a two-tensor
$K$ is called an initial data set $(M,g,K)$ if $g$ and $K$ satisfy
the constraint equations

\begin{eqnarray}
  R_g-|K|_g^2+(tr_g(K))^2=16\pi\rho\nonumber\\
  div_g(K)-d(tr_g(K))=8\pi J
\end{eqnarray}
where $R_g$ is the scalar curvature of the metric $g$, $tr_g(K)$
denotes $g^{ij}K_{ij}$, $\rho$ is the observed energy density, and
$J$ is the observed momentum density.
\definition
      Let $q\in(\frac{1}{2},1]$. We say $(M,g,K)$ is asymptotically
      flat (AF) if it is a initial data set, and there is a compact subset $\widetilde{K}\subset M$ such that $M\setminus \widetilde{K}$ is
      diffeomorphic to $R^3\setminus B_1(0)$ and there exists
      coordinate $\{x^i\}$ such that

      \begin{equation}
        g_{ij}(x)=\delta_{ij}+h_{ij}(x)
      \end{equation}

      \begin{eqnarray}
        h_{ij}(x)=O_5(|x|^{-q})  && K_{ij}(x)=O_1(|x|^{-1-q})
      \end{eqnarray}
Also, $\rho$ and $J$ satisfy

\begin{eqnarray}
  \rho(x)=O(|x|^{-2-2q})&&J(x)=O(|x|^{-2-2q})
\end{eqnarray}

Here, $f=O_k(|x|^{-q})$ means $
  \partial^lf=O(|x|^{-l-q})$
for $l=0,\cdots, k$. $M\setminus \widetilde{K}$ is called an end of
this asymptotically flat manifold.

We can define mass for the asymptotically flat manifolds as follows:
\begin{equation}
    m=\lim_{r\rightarrow\infty}\frac{1}{16\pi}\int_{|x|=r} (h_{ij,j}-h_{jj,i})v_g^id\mu_g
\end{equation}
where $v_g$ and $d\mu_g$ are the normal vector and volume form with
respect to the metric $g$. From \cite{B},we know the mass is well
defined when $q>1/2$.

\definition We say $(M,g,K)$is
asymptotically flat satisfying the Regge-Teitelboim condition
(AF-RT) if it is AF, and $g,K$ satisfy these asymptotically even/odd
conditions
\begin{eqnarray}
  h^{odd}_{ij}(x)=O_2(|x|^{-1-q})&&K_{ij}^{even}(x)=O_1(|x|^{-2-q})
\end{eqnarray}
Also, $\rho$ and $J$ satisfy

\begin{eqnarray}
  \rho^{odd}(x)=O(|x|^{-3-2q})&& J^{odd}(x)=O(|x|^{-3-2q})
\end{eqnarray} where $f^{odd}(x)=f(x)-f(-x)$ and
$f^{even}(x)=f(x)+f(-x)$.

For (AF-RT) manifolds, the center of mass $C$ is defined
    by
    \begin{equation}
      C^\alpha=\frac{1}{16\pi
      m}\lim_{r\rightarrow\infty}(\int_{|x|=r}x^\alpha(h_{ij,i}-h_{ii,j})v^j_gd\mu_g-\int_{|x|=r}(h_{i\alpha}v^i_g-h_{ii}v_g^\alpha)d\mu_g).
    \end{equation}
From \cite{H1}, we know it is well defined.

 The constant mean curvature surface is stable means the second
variation operator has non-negative eigenvalues when restricted to
the functions with $0$ mean value, i.e.
\begin{equation}
  \int_{\Sigma}(|A|^2+Ric(v_g,v_g))f^2d\mu\leq
\int_{\Sigma}|\nabla f|^2 d\mu
\end{equation}
for function $f$  with $\int_{\Sigma}f d\mu=0$, where $A$ is the
second fundamental form, and $Ric(v_g,v_g)$ is the Ricci curvature
in the normal direction with respect to the metric $g$.

We discuss the existence and uniqueness of constant mean curvature
spheres that separate the origin from the infinity in the AF-RT
manifolds. The following two theorems are due to Lan-Hsuan Huang
\cite{H}:

  \begin{theorem}(Existence)\label{29}
  If $(M,g,K)$ is the AF-RT with $q\in (\frac{1}{2},1]$, there
  exists a foliation by spheres $\{\Sigma_R\}$ with constant mean
  curvature $H(\Sigma_R)=\frac{2}{R}+O(R^{-1-q})$ in the
  exterior region of $M$. Each leaf $\Sigma_R$ is a
  $c_0R^{1-q}$-graph over $S_R(C)$ and is strictly stable.
\end{theorem}

Set   $r(x)=(\Sigma (x_i)^2)^{1/2}$.  For the constant mean
curvature sphere $\Sigma$ which separates infinity from $K$, we
define
\begin{eqnarray}
  &&r_0(\Sigma) = \inf\{r(x)|x\in \Sigma\} \nonumber\\
  &&r_1(\Sigma) = \sup\{r(x)|x\in \Sigma\}
\end{eqnarray}

\begin{theorem}\label{33}(Uniqueness)   Assume that $(M,g,K)$ is AF-RT with
$q\in(\frac{1}{2},1]$ and $m>0$. There exists $\sigma_1$ and $C_1$
so that if $\Sigma$ has the following properties:

\begin{itemize}
  \item $\Sigma$ is topologically a sphere
  \item $\Sigma$ has constant mean curvature $H=H(\Sigma_R)$ for
  some $R \geq \sigma_1$
  \item $\Sigma$ is stable
  \item $r_1\leq C_1r_0^{\frac{1}{a}}$ for some $a$ satisfying $\frac{5-q}{2(2+q)}<a\leq1$
\end{itemize}
then $\Sigma=\Sigma_R$.
\end{theorem}

Our main uniqueness result is

\theorem \label{26}Suppose $(M,g,K)$ is AF-RT 3-manifold with
positive mass, and $g$ can be expressed on the end $M\setminus
\widetilde{K}$ as follows:
\begin{equation}
  g_{ij}=\delta_{ij}+h^1_{ij}(\theta)/r+Q
\end{equation}
where $\theta=(\theta_1,\theta_2)$ is the coordinate on $S^2\subset
R^3$. If $g$ satisfies the following properties:

\begin{itemize}
  \item $h^1_{ij}(\theta)\in C^5(S^2)$
  \item $Q=O_5(|x|^{-2})$

\end{itemize}

 Then for any $k>2$, there exists some $\varepsilon>0$ depending on $k$ such that if

 \begin{equation}
   \|h_{ij}(\theta)-\delta_{ij}(\theta)\|_{W^{k,2}(S^2)}\leq \varepsilon ,
 \end{equation}
 there is a compact domain $\widetilde{K}$
such that if a foliation $\{\Sigma\}$ of stable constant mean
curvature spheres
 which separates infinity
from $\widetilde{K}$ have

\begin{equation}
 \lim_{r_0\rightarrow \infty}\frac{\log(r_1(\Sigma))}{r_{0}(\Sigma)^{1/4}}=0
\end{equation}
then this foliation is the same one as in Theorem\ref{29}.

\begin{remark}
  If we replace $\|h_{ij}(\theta)-\delta_{ij}(\theta)\|_{W^{k,2}}\leq
  \varepsilon$ by $\|h_{ij}(\theta)-C\delta_{ij}(\theta)\|_{W^{k,2}}\leq
  \varepsilon$ for any constant $C>0$, we can also get this theorem, but $\varepsilon$ will depend on $k$ and $C$.
\end{remark}

\begin{remark}
  RT condition is needed to apply the theorems of Huang and if we assume the scalar curvature satisfies $R=O(r^{-3-\varepsilon})$ for some $\varepsilon>0$, then we
  do not need the constraint equation.
\end{remark}
\begin{remark}
  Here I can only deal with the case when $q=1$. When $q\in (1/2,1)$
  it seems that $\|h_{ij}(\theta)-\delta_{ij}(\theta)\|_{W^{k,2}(S^2)}\leq
  \varepsilon$ is not a proper condition.
\end{remark}
 The above theorem is about the uniqueness of the foliation. For
the uniqueness of a single CMC sphere we have:

\corollary \label{34}We assume the same condition on the metric as
the above Theorem. Then for any constants $C>0$ and $\beta>0$, there
exist some compact set $K(C,\beta)\subset M$, such that any stable
sphere $\Sigma$ that separates $K(C,\beta)$ from the infinity with
\begin{equation}
 \frac{(\log(r_1(\Sigma)))^{1+\beta}}{r_{0}(\Sigma)^{1/4}}\leq C
\end{equation} belongs to the foliation in Theorem \ref{29}.

The paper is organized much like \cite{QT}: In Section 2 we do
apriori estimate on the stable constant mean curvature sphere based
on the Simon's identity. In Section 3, we introduce blow-down
analysis in three different scales. In Section 4 we recall the
asymptotic analysis from \cite{QT1} and prove a technical lemma. In
Section 5 we introduce the asymptotically harmonic coordinate. In
Section 6 we introduce a sense of the center of mass and prove the
theorem.

\section{Curvature estimates}

From now on let $\Sigma$ be a constant mean curvature sphere in the
asymptotically flat end $(M,g)$which separates the origin from the
infinity. First we have the following estimate as Lemma 5.2 in
\cite{HY}.

\lemma Let $X=x^i\frac{\partial}{\partial x^i}$ be the Euclidean
coordinate vectorfield and $r=(\Sigma (x^i)^2)^{1/2}$ and with
respect to the metric $g$, $v$ is the outward normal vector field ,
$d\mu$ is the volume form of $\Sigma$. Then we have the estimate:
\begin{equation}
  \int_\Sigma<X,v>^2r^{-4}d\mu\leq H^2|\Sigma|
\end{equation}

Moreover for each $a\geq a_0>2$ and $r_0$ sufficiently large , we
have:

\begin{equation}
  \int_\Sigma r^{-a}d\mu\leq C(a_0)r_0^{2-a}H^2|\Sigma|
\end{equation}

Proof. Because the mean curvature $H$ is constant, then for some
smooth vector field $Y$ on $\Sigma$ , we have the divergence
formula:
\begin{equation}
  \int_\Sigma div_\Sigma Y d\mu=H\int_\Sigma<Y,v>d\mu.
\end{equation}
We choose $Y=Xr^{-a}$ , $a\geq2$ and $e_\alpha$ is the orthonormal
basis on $\Sigma$ , $\alpha=1,2$. Suppose
$e_\alpha=a_\alpha^i\frac{\partial}{\partial x^i}$, it is obvious
that $a_\alpha^i$ is bounded because the manifold is asymptotically
flat. Then we have:

\begin{eqnarray}
  div_\Sigma
  Y&&=div_\Sigma(Xr^{-a})=<\nabla_{e_\alpha}(Xr^{-a}),e_{\alpha}>\nonumber\\
  &&=r^{-a}div_\Sigma X-ar^{-a-2}a_\alpha^ia_\alpha^jx^ix^j+O(r^{-a-q})\nonumber\\
  &&=r^{-a}div_\Sigma
  X-\alpha r^{-a-2}|X^\tau|^2+O(r^{-a-q})
\end{eqnarray}
where $X^\tau$ is the tangent projection of $X$.

\begin{eqnarray}
  |div_\Sigma X-2|=O(r^{-q})
\end{eqnarray}

Note that$|X^\tau|^2=r^2-<X,v>^2+O(r^{2-q})$ , then combine all of
these we have:

\begin{eqnarray}\label{13}
  &&|(2-a)\int_\Sigma
  r^{-a}d\mu+a\int_{\Sigma}<X,v>^2r^{-a-2}d\mu-H\int_{\Sigma}<X,v>r^{-a}d\mu|\nonumber\\
  &&\leq
  C\int_{\Sigma}r^{-a-q}d\mu
\end{eqnarray}

Choosing $a=2$ , from H\"{o}lder inequality , we have:

\begin{equation}\label{12}
  \int_{\Sigma}<X,v>^2r^{-4}d\mu\leq\frac{1}{4}H^2|\Sigma|+C\int_{\Sigma}r^{-2-q}d\mu
\end{equation}

then choose $a=2+q$ ,

\begin{equation}
  \int_{\Sigma}r^{-2-q}d\mu\leq4r_0^{-q}(\int_{\Sigma}<X,v>^2r^{-4}d\mu+H^2|\Sigma|+C\int_{\Sigma}r^{-2-q}d\mu)
\end{equation}

then combine this with (\ref{12}),we have:

\begin{equation}
  \int_{\Sigma}<X,v>^2r^{-4}d\mu\leq H^2|\Sigma|
\end{equation}

then again from (\ref{13}), we have for $a\geq a_0>2$, we derive:

\begin{equation}
  \int_{\Sigma}r^{-a}\leq C(a_0-2)^{-1}r_0^{2-a}H^2|\Sigma|
\end{equation}

 Then we can derive the
integral estimate for $|{\AA}|$ from the stability of the surface as
in \cite{HY} Proposition 5.3, i.e. we have

\lemma Suppose $\Sigma$ is a stable constant mean curvature sphere
in the asymptotically flat manifold. We have for $r_0$ sufficiently
large
\begin{equation}
\int_{\Sigma}|{\AA}|^2d\mu\leq Cr_0^{-q}
\end{equation}
\begin{equation}
  H^2|\Sigma|\leq C
\end{equation}
\begin{equation}
  \int_{\Sigma}H^2d\mu=16\pi+O(r_0^{-q})
\end{equation}

Proof. Since $\Sigma$ is stable , we have

\begin{equation}
  \int_{\Sigma}|\nabla f|^2d\mu\geq
  \int_{\Sigma}(|A|^2+Ric(v,v))f^2d\mu
\end{equation}
for any function $f$ , with$\int_{\Sigma}f d\mu=0$, where $A$ is the
second fundamental form of $\Sigma$ and $Ric$ is the Ricci curvature
of $M$

Choose $\psi$ to be a conformal map of degree 1 from $\Sigma$ to the
standard $S^2$ in $R^3$. Each component $\psi_i$ of $\psi$ can be
chosen such that $\int \psi_id\mu=0$ , see \cite{LY} . We have for
each $\psi_i$

\begin{equation}
  \int_\Sigma|\nabla \psi_i|^2d\mu=\frac{8\pi}{3}
\end{equation}

since $\sum \psi_i^2\equiv1$ we conclude that

\begin{equation}
  \int_{\Sigma}|A|^2+Ric(v,v)d\mu\leq 8\pi
\end{equation}

From Gauss equation

\begin{equation}\label{14}
  \frac{1}{2}|A|^2+Ric(v,v)-\frac{1}{2}R+K=\frac{1}{2}H^2
\end{equation}
 we have:

\begin{equation}
  |A|^2+Ric(v,v)=\frac{1}{2}|{\AA}|^2+\frac{3}{4}H^2+\frac{1}{2}R-K
\end{equation}
where $K$ is the Gauss curvature of $\Sigma$ and ${\AA}$ is defined
as ${\AA}_{ij}=A_{ij}-\frac{H}{2}g_{ij}$

Then we have:

\begin{equation}
  \int_{\Sigma}\frac{1}{2}|{\AA}|^2+\frac{3}{4}H^2|\Sigma|\leq 12\pi
  +r_0^{-q}H^2|\Sigma|
\end{equation}
because $R=O(r^{-2-2q})$.

So we have $H^2|\Sigma|\leq 16\pi$.

Using the Gauss equation in a different way, we have

\begin{eqnarray}
  &&\int_{\Sigma}|{\AA}|^2d\mu=\int_{\Sigma}|A|^2-\frac{H^2}{2}d\mu\nonumber\\
&&=\frac{1}{2}\int_{\Sigma}|A|^2+Ric(v,v)d\mu+\frac{1}{2}\int_{\Sigma}R-3Ric(v,v)-2Kd\mu\nonumber\\
&&\leq\int_{\Sigma}r^{-2-q}d\mu\nonumber\\
&&=O(r_0^{-q}).
\end{eqnarray}

Then from Gauss equation (\ref{14}) again, we have:

\begin{equation}
\int_{\Sigma}H^2d\mu=4\int_{\Sigma}Kd\mu+O(r_0^{-q})=16\pi+O(r_0^{-q})
\end{equation}

 \lemma Suppose that M is a constant mean curvature surface in an
asymptotically flat end $(R^3\setminus B_1(0),g)$. Then

\begin{equation}
  \int_{\Sigma}H_e^2d\mu_e=16\pi+O(r_0^{-q})
\end{equation}

\proof We follow the calculation of Huisken and Ilmanen \cite{HI},
\begin{equation}\label{1}
    g_{ij}=\delta_{ij}+h_{ij}
\end{equation}

Suppose
\begin{equation}   g_{ij}|_\Sigma=f_{ij}    ,
\delta_{ij}|_\Sigma=\varepsilon_{ij}
\end{equation}
$f^{ij}$  and
  $\epsilon^{ij}$  are the corresponding inverse matrices. $v,\omega,A,H,d\mu$
represents the normal vector , the dual form of $v$, the second
fundamental form , the mean curvature and the volume form of
$\Sigma$ in the metric $g$. And
$v_e,\omega_e,A_e,H_e,\mu_e$
represents the corresponding ones in Euclidean metric. Through easy
calculation, we have
\begin{eqnarray}
  &&f^{ij}-\varepsilon^{ij}=-f^{ik}h_{kl}f^{lj}\pm C|h|^2 \label{2}\\
  &&g^{ij}-\delta^{ij}=-g^{ik}h_{kl}g^{lj}\pm C|h|^2 \label{3}
\end{eqnarray}
\begin{eqnarray}
  \omega=\frac{\omega_e}{|\omega_e|} &&
  v^i=g^{ij}\omega_j
\end{eqnarray}
\begin{eqnarray}\label{4}
  (\omega_e)_i=\omega_i\pm C|P| &
v^i_e=v^i+C|h| & 1-|\omega_e|=\frac{1}{2}h_{ij}v^iv^j
\end{eqnarray}
\begin{equation}
  \Gamma_{ij}^k=\frac{1}{2}g^{kl}(\overline{\nabla}_ih_{jl}+\overline{\nabla}_jh_{il}-\overline{\nabla}_lh_{ij})\pm
  C|h|\pm C|\overline{\nabla} h| \label{5}
\end{equation}
and $\Gamma_{ij}^k$is the Christoffel symbol for
$\overline{\nabla}-\overline{\nabla}_e$ ,where we denote the gradient for the metric $g$ and $\delta$ by $\overline{\nabla}$ and $\overline{\nabla}_e$.


We have the formula:
\begin{equation}
  |\omega_e|_gA_{ij}=(A_e)_{ij}-(\omega_e)_k\Gamma_{ij}^k
  \label{6}
\end{equation}

So we have

\begin{eqnarray}
  &&H-H_e=f^{ij}A_{ij}-\varepsilon^{ij}(A_e)_{ij} \nonumber\\
   &&=
   (f^{ij}-\varepsilon^{ij})A_{ij}+\varepsilon^{ij}A_{ij}(1-|\omega_e|_g)+\varepsilon^{ij}(|\omega_e|_gA_{ij}-(A_e)_{ij})
\end{eqnarray}
from (\ref{2})(\ref{3})(\ref{4}), we have
\begin{equation}
  \varepsilon^{ij}A_{ij}(1-|\omega_e|_g)=\frac{1}{2}Hv^iv^jh_{ij}\pm
  C|h|^2|A|
\end{equation}
and using (\ref{2})(\ref{3})(\ref{4})(\ref{5})(\ref{6})we have:

\begin{eqnarray}
  &&  \varepsilon^{ij}(|\omega_e|A_{ij}-(A_e)_{ij})\nonumber\\
  &&=-\varepsilon^{ij}(\omega_e)_k\Gamma_{ij}^k\nonumber \\
    &&=-\frac{1}{2}f^{ij}\omega_kg^{kl}(\overline{\nabla}_ih_{jl}+\overline{\nabla}_jh_{il}-\overline{\nabla}_lh_{ij})\pm C|h||\overline{\nabla} h| \nonumber\\
    &&=-f^{ij}v^l\overline{\nabla}_ih_{jl}+\frac{1}{2}f^{ij}v^l\overline{\nabla}_lh_{ij}\pm
C|h||\overline{\nabla} h|
\end{eqnarray}

At last , we have
\begin{eqnarray}\label{16}
  &&H-H_e=-f^{ik}h_{kl}f^{lj}A_{ij}+\frac{1}{2}Hv^iv^jh_{ij}-f^{ij}v^l\overline{\nabla}_ih_{jl}\nonumber\\
  &&+\frac{1}{2}f^{ij}v^l\overline{\nabla}_lh_{ij}\pm
  C|h||\overline{\nabla} h|\pm C|h|^2|A|
\end{eqnarray}

\begin{eqnarray}
&&\int_{\Sigma}H_e^2d\mu_e= (1+O(r_0^{-q}))\int_{\Sigma}H_e^2d\mu\nonumber\\
&&\leq(1+O(r_0^{-q}))(\int_{\Sigma}H^2d\mu+\int_{\Sigma}(H_e-H)^2+2|H(H_e-H)|d\mu)\nonumber\\
&&\leq(1+O(r_0^{-q}))(16\pi+O(r_0^{-q})+\int_{\Sigma}(H_e-H)^2\nonumber\\
&&+(\int_{\Sigma}H^2d\mu)^{\frac{1}{2}}(\int_{\Sigma}(H_e-H)^2d\mu)^{\frac{1}{2}})
\end{eqnarray}

\begin{eqnarray}
  &&\int (H_e-H)^2d\mu\leq\int
  O(|x|^{-2q})|A|^2+H^2O(|x|^{-2q})+O(|x|^{-2-2q})d\mu\nonumber\\
  &&\leq\int O(|x|^{-2q})H^2+ O(|x|^{-2q})|{\AA}|^2+
  O(|x|^{-2-2q})d\mu\nonumber\\
  &&=O(r_0^{-2q})
\end{eqnarray}
so we have

\begin{equation}
  \int_{\Sigma}H_e^2d\mu_e\leq16\pi+O(r_0^{-q})
\end{equation}

On the other hand, by Euler formula,
\begin{equation}
  K_e=\frac{1}{4}H_e^2-\frac{1}{2}|{\AA}_e|^2.
\end{equation}

So we have
\begin{equation}
  \int H_e^2d\mu_e\geq16 \pi
\end{equation}
which implies:
\begin{equation}
 \int_{\Sigma} H_e^2d\mu_e=16\pi+O(r_0^{-q})
\end{equation}

Based on Michael and Simon, we have the following Sobolev
inequality.

\lemma Suppose that $\Sigma$ is a constant mean curvature surface in
an asymptotically flat end$(R^3\setminus B_1(0),g)$ with
$r_0(\Sigma)$ sufficiently large, and that $\int_\Sigma H^2\leq C$.
Then
\begin{equation}\label{15}
  (\int_\Sigma f^2d\mu)^{\frac{1}{2}}\leq C(\int_\Sigma|\nabla
  f|d\mu+\int_
  \Sigma H|f|d\mu).
\end{equation}
\proof Note that it is valid for the surface in Euclidean Space. So
by the uniform equivalence of the metric $g$ and $\delta$ , we have:
\begin{eqnarray}
  (\int|f|^2d\mu)^{\frac{1}{2}}\leq
  C(\int|f|^2d\mu_e)^{\frac{1}{2}}\leq C(\int |\nabla f|+H|f|+|H-H_e||f|d\mu)
\end{eqnarray}
To bound the last term on the right , we have:
\begin{eqnarray}
  \int|H-H_e||f|d\mu&\leq& \int
  O(|x|^{-q})|A||f|+O(|x|^{-q})H|f|\nonumber\\
  &&+O(|x|^{-1-q})|f|d\mu\nonumber\\
  &\leq& O(r_0^{-q})\int
  H|f|+(\int|{\AA}|^2d\mu)^{\frac{1}{2}}O(r_0^{-q})\|f\|_{L^2}\nonumber\\
  &&+O(r_0^{-q})\|f\|_{L^2}
\end{eqnarray}
So we can choose $r_0$ sufficiently large and get the desired
result.

\lemma Suppose that $\Sigma$ is a constant mean curvature surfaces
in an asymptotically flat end $(R^3\setminus B_1(0),g)$ with
$r_0(\Sigma)$ sufficiently large, then:

\begin{equation}
  C_1H^{-1}\leq diam(\Sigma) \leq C_2H^{-1}
\end{equation}

In particular, if the surface $\Sigma$ separates the infinity from the compact part, then:

\begin{equation}
   C_1H^{-1}\leq r_1(\Sigma) \leq C_2H^{-1}
\end{equation}
Proof.  We already know that:

\begin{equation}
 \int_{\Sigma}H_e^2d\mu_e=16\pi+O(r_0^{-q})
\end{equation}

Then from \cite{L} Lemma 1.1, we know that

\begin{equation}
  \sqrt{\frac{2|\Sigma|_e}{F(\Sigma)}}\leq diam(\Sigma)\leq C\sqrt{|\Sigma|_eF(\Sigma)}
\end{equation}
where $F(\Sigma)=\frac{1}{2}\int_{\Sigma}H_e^2$ is the Willmore
functional and $|\Sigma|_e$ is the volume of $\Sigma$ with respect
to the Euclidean metric. But the Euclidean metric is uniformly
equivalent to $g$, so we get the result.

Now to get the pointwise estimate for ${\AA}$ ,we use the Simons
identity and the Moser's iteration argument.

\lemma(Simons identity \cite{SSY}) Suppose $N$ is a hypersurface in
a Riemannian manifold $(M,g)$ , then the second fundamental form
satisfies the following identity:

\begin{eqnarray}
 && \Delta
  A_{ij}=\nabla_i\nabla_jH+HA_{ik}A_{jk}-|A|^2A_{ij}+HR_{3i3j}-A_{ij}R_{3k3k}
  +A_{jk}R_{klil}\nonumber\\
  &&+A_{ik}R_{kljl}-2A_{lk}R_{iljk}+\overline{\nabla}_jR_{3kik}+\overline{\nabla}_kR_{3ijk}
\end{eqnarray}
where $R_{ijkl}$ and $\overline{\nabla}$ are the curvature and
gradient operator of $(M,g)$, then from this we easily deduce for
constant mean curvature surface we have the next inequality for
${\AA}$ :

\begin{eqnarray}
  &&-|{\AA}|\Delta|{\AA}|\leq
  |{\AA}|^4+CH|{\AA}|^3+CH^2|{\AA}|^2+C|{\AA}|^2|x|^{-2-q}\nonumber\\
  &&+CH|{\AA}||x|^{-2-q}+C|{\AA}||x|^{-3-q}
\end{eqnarray}



We also need an inequality for $\nabla\AA$ because we also want to
estimate the higher derivative:

\begin{eqnarray}
  &&-|\nabla{\AA}|\Delta|\nabla{\AA}|\leq
  C|\nabla{\AA}|^2(|\AA|^2+H|\AA|+H^2+O(|x|^{-2-q}))\\
  &&+|\nabla{\AA}|((|{\AA}|^2+H|\AA|+H^2)O(|x|^{-2-q})+(|{\AA}|+H)O(|x|^{-3-q})+O(|x|^{-4-q}))\nonumber
\end{eqnarray}

\lemma
\begin{equation}
  \|{\AA}^2\|_{L^2}+\|\nabla|{\AA}|\|_{L^2}+\|\nabla {\AA}\|_{L^2}+\|H|{\AA}|\|_{L^2}\leq Cr_0^{-1-q}
\end{equation}
Proof. See \cite{H} Lemma 4.5

Then we can get the pointwise estimates for ${\AA}$ and $\nabla{\AA}$ .

\theorem \cite{QT}\label{21}Suppose that $(R^3\setminus B_1(0),g)$
is an asymptotically flat end. Then there exist positive numbers
$\sigma_0$, $\delta_0$ such that for any constant mean curvature
surface in the end, which separates the infinity from the compact
part, we have:

\begin{equation}
  |{\AA}|^2(x)\leq C|x|^{-2}\int_{B_{\delta_0|x|}(x)}|{\AA}|^2d\mu+C|x|^{-2-2q}\leq C|x|^{-2}r_0^{-q}
\end{equation}

\begin{equation}
  |\nabla{\AA}|^2(x)\leq C|x|^{-2}\int_{B_{\delta_0|x|}(x)}|\nabla{\AA}|^2d\mu+C|x|^{-4-2q}\leq C|x|^{-2}r_0^{-2-2q}
\end{equation}
provided that $r_0\geq \sigma_0$.

Proof. In the Sobolev inequality (\ref{15}) we take $f=u^2$ , then we get:

\begin{eqnarray}
  &&(\int_{\Sigma}u^4d\mu)^{\frac{1}{2}}\leq C(2\int_{\Sigma}|u||\nabla u|d\mu+\int_{\Sigma }Hu^2d\mu)\nonumber\\
  &&\leq C(\int_{\Sigma}u^2)^{\frac{1}{2}}(\int_{\Sigma}|\nabla u|^2d\mu)^{\frac{1}{2}}+C(\int_{supp(u)}H^2d\mu)^{\frac{1}{2}}(\int_{\Sigma}u^4d\mu)^{\frac{1}{2}}
\end{eqnarray}

\lemma For any $\varepsilon>0$, we can find a uniform $\delta_0$ sufficiently small such that if for any $x\in\Sigma$ , we have that:

\begin{equation}
  \int_{B_{\delta_0|x|}(x)}H^2\leq\varepsilon
\end{equation}
Proof. In fact we need only to prove that there exist $C$

\begin{equation}
  |B_{\delta_0|x|}(x)|\leq C\delta_0^2|x|^2
\end{equation}
because then,

\begin{equation}
  H^2|B_{\delta_0|x|}(x)|\leq C\delta_0^2|x|^2H^2\leq C\delta_0^2
\end{equation}

From \cite{L} the proof of lemma 1.1, we know that, for any $x\in
\Sigma$, $B_\sigma(x)$ denotes the Euclidean ball of radius $\sigma$
with center $x$ in $R^3$, $\Sigma_\sigma=\Sigma\cap B_\sigma(x)$,
then there exists $C$ such that for $0<\sigma\leq\rho<\infty$

\begin{equation}
  \sigma^{-2}|\Sigma_\sigma|\leq C(\rho^{-2}|\Sigma_\rho|+F(\Sigma_\rho))
\end{equation}
where $F(\Sigma_\rho)$ is the Willmore functional. $C$ doesn't depend on $\Sigma , \sigma , \rho$.

Let $\rho\rightarrow\infty$ , $\rho^{-2}|\Sigma_\rho|\rightarrow0$,
so we have:

\begin{equation}
  \sigma^{-2}|\Sigma_\sigma|\leq CF(\Sigma)\leq C
\end{equation}
so we prove the lemma.

So if $supp(u)\subset B_{\delta_0|x|}(x)$, we have the following
scaling invariant Sobolev inequality:

\begin{equation}
  (\int_{\Sigma}u^4d\mu)^{\frac{1}{2}}\leq C(\int_{\Sigma}u^2)^{\frac{1}{2}}(\int_{\Sigma}|\nabla u|^2d\mu)^{\frac{1}{2}}
\end{equation}

\lemma \cite{QT} Suppose that a nonnegative function $v\in L^2$
solves

\begin{equation}
  -\Delta v\leq fv+h
\end{equation}
on $B_{2R}(x_0)$, where

\begin{equation}
  \int_{B_{2R}(x_0)}f^2d\mu \leq CR^{-2}
\end{equation}
and $h\in L^2(B_{2R}(x_0))$. And suppose that

\begin{equation}
  (\int_{\Sigma}u^4d\mu)^{\frac{1}{2}}\leq C(\int_{\Sigma}u^2)^{\frac{1}{2}}(\int_{\Sigma}|\nabla u|^2d\mu)^{\frac{1}{2}}
\end{equation}
holds for all $u$ with support inside $B_{2R}(x_0)$. Then

\begin{equation}
  \sup_{B_R(x_0)}v\leq CR^{-1}\|v\|_{L^2(B_{2R}(x_0))}+CR\|h\|_{L^2(B_{2R(x_0)})}
\end{equation}

See \cite{QT} Lemma 2.6 for the proof of this lemma.

Then we find that:

\begin{eqnarray}
  &&-\Delta|{\AA}|\leq(|{\AA}|^2+H^2+H|{\AA}|+C|x|^{-2-q})|{\AA}|+CH|x|^{-2-q}+C|x|^{-3-q}\nonumber\\
  &&=f_1|\AA|+h_1
\end{eqnarray}

\begin{eqnarray}
  &&-\Delta|\nabla{\AA}|\leq
  C|\nabla{\AA}|(|\AA|^2+H|\AA|+H^2+O(|x|^{-3}))\nonumber\\
  &&+((|{\AA}|^2+H|\AA|+H^2)O(|x|^{-3})+(|{\AA}|+H)O(|x|^{-4})+O(|x|^{-5}))\nonumber\\
  &&=f_2|\nabla{\AA}|+h_2.
\end{eqnarray}
We need to prove that
$\|f_1\|^2_{L^2(B_{2\delta_0|x|}(x))},\|f_2\|^2_{L^2(B_{2\delta_0|x|}(x))}\leq
C|x|^{-2}$ , see \cite{QT} Theorem 2.5 for the proof. and it is easy
to show that $\|h_1\|^2_{L^2(B_{2\delta_0|x|}(x))}=O(|x|^{-4-2q})$
and $\|h_2\|^2_{L^2(B_{2\delta_0|x|}(x))}=O(|x|^{-6-2q})$.

\begin{remark}
  We can also do the same kind of estimate for $\nabla^2{\AA}$,
  where we need the third derivative of curvature. It is needed by the $C^{2,\alpha}$ convergence of the surface in the next section. This is the
  reason why we require the metric $g$ to be smooth up to 5th order.
\end{remark}

\section{Blow down analysis}

Now like \cite{QT}, we blow down the surface in three different
scales. First we consider

\begin{equation}
  \widetilde{N}=\frac{1}{2}HN=\{\frac{1}{2}Hx:x\in N\}
\end{equation}

Suppose that there is a sequence of constant mean curvature surfaces
$\{N_i\}$ such that

\begin{equation}
  \lim_{i\rightarrow\infty}r_0(N_i)=\infty
\end{equation}
we have known that

\begin{equation}
   \lim_{i\rightarrow\infty}\int_{N_i}H_e^2d\sigma=16\pi
\end{equation}

Hence, by the curvature estimates established in the previous
section combining the proof of Theorem 3.1 in \cite{L}, we have

\lemma \label{18}Suppose that $\{N_i\}$ is a sequence of constant
mean curvature surfaces in a given asymptotically flat end
$(R^3\setminus B_1(0),g)$ and that

\begin{equation}
  \lim_{i\rightarrow\infty}r_0(N_i)=\infty.
\end{equation}
And suppose that $N_i$ separates the infinity from the compact part.
Then, there is a subsequence of $\{\widetilde{N}_i\}$ which
converges in Gromov-Hausdorff distance to a round sphere $S_1^2(a)$
of radius $1$ and centered at $a\in R^3$. Moreover,the convergence
is in $C^{2,\alpha}$ sense away from the origin.

Then, we use a smaller scale $r_0$ to blow down the surface

\begin{equation}
  \widehat{N}=r_0(N)^{-1}N=\{r_0^{-1}x:x\in N\}.
\end{equation}

\lemma \label{17}Suppose that $\{N_i\}$ is a sequence of constant
mean curvature surfaces in a given asymptotically flat end
$(R^3\setminus B_1(0),g)$ and that

\begin{equation}
   \lim_{i\rightarrow\infty}r_0(N_i)=\infty.
\end{equation}

And suppose that

\begin{equation}
  \lim_{i\rightarrow\infty}r_0(N_i)H(N_i)=0.
\end{equation}

Then there is a subsuquence of $\{\widehat{N}_i\}$ converges to a
2-plane at distance $1$ from the origin. Moreover the convergence is
in $C^{2,\alpha}$ in any compact set of $R^3$.

We must understand the behavior of the surfaces $N_i$ in the scales
between $r_0(N_i)$ and $H^{-1}(N_i)$. We consider the scale $r_i$
such that

\begin{eqnarray}
  \lim_{i\rightarrow\infty}\frac{r_0(N_i)}{r_i}=0 &&
  \lim_{i\rightarrow\infty}r_iH(N_i)=0
\end{eqnarray}

and blow down the surfaces

\begin{equation}
  \overline{N}_i=r_i^{-1}N=\{r_i^{-1}x:x\in N\}.
\end{equation}

\lemma \label{19}Suppose that $\{N_i\}$ is a sequence of constant
mean curvature surfaces in a given asymptotically flat end
$(R^3\setminus B_1(0),g)$ and that

\begin{equation}
  \lim_{i\rightarrow\infty}r_0(N_i)=\infty
\end{equation}

And suppose that $r_i$ are such that

\begin{eqnarray}
  \lim_{i\rightarrow\infty}\frac{r_0(N_i)}{r_i}=0 &&
  \lim_{i\rightarrow\infty}r_iH(N_i)=0
\end{eqnarray}

Then there is a subsequence of $\{\overline{N}_i\}$ converges to a
2-plane at the origin in Gromov-Hausdorff distance. Moreover the
convergence is $C^{2,\alpha}$ in any compact subset away from the
origin.

\section{Asymptotically analysis}
First we revise  Proposition 2.1 in \cite{QT1}. We prove a different
version. Let us denote:

\begin{equation}
  \|u\|^2_{1,i}=\int_{[(i-1)L,iL]\times S^1}|u|^2+|\nabla u|^2dtd\theta
\end{equation}

\lemma \label{20}Suppose $u\in W^{1,2}(\Sigma,R^k)$ satisfies

\begin{equation}
  \Delta u+ A\cdot\nabla u+B\cdot u=h
\end{equation}
in $\Sigma$, where $\Sigma=[0,3L]\times S^1$. And suppose that $L$
is given and large. Then there exists a positive number $\delta_0$
such that if

\begin{equation}
  |h|_{L^2(\Sigma)}\leq \delta_0 \max_{1\leq i\leq3}|u|_{1,i}
\end{equation}
and

\begin{eqnarray}
  |A|_{L^\infty(\Sigma)}\leq \delta_0&&|B|_{L^\infty(\Sigma)}\leq \delta_0
\end{eqnarray}
then,

(a)$\|u\|_{1,3}\leq e^{-\frac{1}{2}L}\|u\|_{1,2}$ implies
$\|u\|_{1,2}< e^{-\frac{1}{2}L}\|u\|_{1,1}$

(b)$\|u\|_{1,1}\leq e^{-\frac{1}{2}L}\|u\|_{1,2}$ implies
$\|u\|_{1,2}< e^{-\frac{1}{2}L}\|u\|_{1,3}$

(c)If both $\int_{L\times S^1}ud\theta$ and $\int_{2L\times
S^1}ud\theta\leq \delta_0\max_{1\leq i\leq3}\|u\|_{1,i}$, then
either $\|u\|_{1,2}< e^{-\frac{1}{2}L}\|u\|_{1,1}$ or $\|u\|_{1,2}<
e^{-\frac{1}{2}L}\|u\|_{1,3}$

Proof. Suppose that $u\in W^{1,2}(\Sigma)$ and $u$ is harmonic, we
can deduce that $u$ satisfies (a)(b)(c')with

(c')If both $\int_{L\times S^1}ud\theta$ and $\int_{2L\times
S^1}ud\theta=0$, then either $\|u\|_{1,2}<
e^{-\frac{1}{2}L}\|u\|_{1,1}$ or $\|u\|_{1,2}<
e^{-\frac{1}{2}L}\|u\|_{1,3}$

A harmonic function $u$ can be written as:

\begin{equation}
  u=a_0+b_0t+\sum_{n=1}^\infty\{e^{nt}(a_n\cos n\theta+b_n\sin n\theta)+e^{-nt}(a_{-n}\cos n\theta+b_{-n}\sin n\theta)\}
\end{equation}

Then it follows that:

\begin{eqnarray}
  &&\|u\|_{1,i}^2=2\pi((a_0^2+b_0^2)L+a_0b_0L^2(2i-1)+\frac{1}{3}b_0^2L^3(3i^2-3i+1))\nonumber\\
  &&+\frac{\pi}{2}\sum_{n=1}^\infty\{\frac{e^{2nL-1}}{n}(e^{2(i-1)nL}(a_n^2+b_n^2)+e^{-2niL}(a_{-n}^2+b_{-n}^2))
  +4L(a_na_{-n}+b_nb_{-n})\}\nonumber\\
  &&+\pi\sum_{n=1}^\infty\{\frac{e^{2nL-1}}{n}(e^{2(i-1)nL}(n^2a_n^2+n^2b_n^2)+e^{-2niL}(n^2a_{-n}^2+n^2b_{-n}^2))\nonumber
  \\
  &&+4L(n^2a_na_{-n}+n^2b_nb_{-n})\}
\end{eqnarray}
$i=1,2,3$

If $L$ is fixed and sufficiently large, then we have

\begin{equation}
  \|u\|_{1,2}^2<\frac{1}{2}(e^L\|u\|^2_{1,3}+e^{-L}\|u\|_{1,1}^2)
\end{equation}
which implies (a). We get (b) in the same way. For (c'), we have
$a_0=b_0=0$ then we have

\begin{equation}
  \|u\|_{1,2}^2<\frac{1}{2}e^{-L}(\|u\|^2_{1,3}+\|u\|_{1,1}^2)
\end{equation}
which implies (c')

The second step is to pass limits. If the proposition were false,
then one would have a sequence of $\delta_k\rightarrow0$ and a
sequence of solution $u_k$ with $\|h_k\|_{L^2}\leq\delta_k$
$|A_k|\leq\delta_k$ and $|B_k|\leq\delta_k$ solves:

\begin{eqnarray}
  \Delta u_k+ A_k\cdot\nabla u_k+B_k\cdot u_k=h_k
\end{eqnarray} We may assume $\max_{1\leq i\leq3}\|u_k\|_{1,i}=1$ otherwise we
can normalize them. Then we know that there is a subsequence that
converges to some $u\in W^{1,2}(\Sigma)$ weakly. And $u$ is a
harmonic function. From the interior $W^{2,p}$ estimate we know the
convergence is strongly $W^{1,2}$ in $I_2$, which implies that $u$
is not trivially zero. Because, with the assumption of the proof by
contradiction, the middle one is the largest.

And because $u_i\rightharpoonup u$ weakly in $W^{1,2}(\Sigma)$
sense. So $u_i\rightharpoonup u$ in $W^{1,2}(I_1)$ and
$W^{1,2}(I_3)$ sense, then we have:

\begin{equation}
  \liminf_{i\rightarrow\infty}\|u_i\|_{1,1}\geq \|u\|_{1,1},\liminf_{i\rightarrow\infty}\|u_i\|_{1,3}\geq \|u\|_{1,3}
\end{equation}
and

\begin{equation}
  \lim_{i\rightarrow\infty}\|u_i\|_{1,2}=\|u\|_{1,2}
\end{equation}
then $u_i$ converges to some non-trivial harmonic function $u$ which
violates one of (a)(b) or (c), which proves the lemma.

From now on we assume $q=1$.

Given a surface $N$ in $R^3$, recall from, for example, (8.5) in
\cite{KA}, that

\begin{equation}
  \Delta_e v+|\nabla_e v|^2v=\nabla_e H_e
\end{equation}
where $v$ is the Gauss map from $N\rightarrow S^2$. For the constant
mean curvature surfaces in the asymptotically flat end
$(R^3\setminus B_1(0),g)$, we have

\lemma

\begin{equation}
  |\nabla_e H_e|(x)\leq C|x|^{-2}r_0^{-1}
\end{equation}

\proof Because the metric $g$ and the Euclidean metric are uniformly
equivalent. So we just prove that

\begin{equation}
  |\nabla H_e|(x)\leq C|x|^{-2}r_0^{-1}
\end{equation}

From (\ref{16}), we know that:

\begin{eqnarray}
  &&|\nabla H_e|\leq |\overline{\nabla} h_{ij}||A|+|h_{ij}||A|^2+|h_{ij}||\nabla
  {\AA}_{ij}|+H|A||h_{ij}|+H|\overline{\nabla}h_{ij}|\nonumber\\
  &&+|A||\overline{\nabla}h_{ij}|+|\overline{\nabla}^2h|\nonumber\\
  &&\leq |x|^{-2}r_0^{-1}
\end{eqnarray}

Suppose $\Sigma$ is a constant mean curvature surface in the
asymptotically flat end. Set

\begin{equation}
  A_{r_1,r_2}=\{x\in\Sigma:r_1\leq|x|\leq r_2\}
\end{equation}
and $A^0_{r_1,r_2}$ stands for the standard annulus in $R^2$. We are
concerned with the behavior of $v$ on
$A_{Kr_0(\Sigma),sH^{-1}(\Sigma)}$ of $\Sigma$ where $K$ will be
fixed large and $s$ will be fixed small. The lemma below gives us a
good coordinate on the surface.

\lemma Suppose $\Sigma$ is a constant mean curvature surface in a
given asymptotically flat end $(R^3\setminus B_1(0),g)$. Then, for
any $\varepsilon>0$ and $L$ fixed and large, there are $M$,$s$ and
$K$ such that, if $r_0\geq M$ and $Kr_0(\Sigma)<r<sH^{-1}(\Sigma)$,
then $(r^{-1}A_{r,e^{L}r},r^{-2}g_e)$ may be represented as
$(A^0_{1,e^{L}},\overline{g})$ and

\begin{equation}
  \|\overline{g}-|dx|^2\|_{C^1(A^0_{1,e^{L}})}\leq \varepsilon.
\end{equation}

In other words, in the cylindrical coordinates $(S^1\times[\log r,
L+\log r,\overline{g}_c])$

\begin{equation}
  \|\overline{g}_c-(dt^2+d\theta^2)\|_{C^1(S^1\times[\log r,L+\log r])}\leq\varepsilon
\end{equation}

Proof. Suppose this is not true. Then we can assume that such $K$
(or such $s$) cannot be found. Then by Lemma \ref{17}. for some
$\varepsilon_0>0$, there is a
 sequence $\Sigma_n$ with
$r_0(\Sigma_n)\rightarrow\infty$, and
$\widetilde{l}_n\rightarrow\infty$ , such that:

\begin{equation}
  ((Kr_0e^{\widetilde{l}_nL})^{-1}A_{Kr_0e^{\widetilde{l}_nL},Kr_0e^{(\widetilde{l}_n+1)L}}, (Kr_0e^{\widetilde{l}_nL})^{-2}g_e)
\end{equation}
is not $\varepsilon_0$ close to $(A^0_{1,e^{L}},\overline{g})$.

By Lemma \ref{18}. We know that

\begin{equation}
  \frac{Kr_0e^{\widetilde{l}_nL}}{sH^{-1}(\Sigma_n)}\rightarrow 0
\end{equation}
must hold because we have choose $s$ sufficiently small.

So if we assume $r_n=Kr_0e^{\widetilde{l}_nL}$, we have:

\begin{equation}
  \lim_{n\rightarrow\infty}\frac{r_n}{Kr_0}=\infty,
  \lim_{n\rightarrow\infty}\frac{r_n}{sH^{-1}}=0
\end{equation}

We blow down the surface using $r_n$, and have a contradiction with
Lemma \ref{19}. This proves the lemma.

Now consider the cylindrical coordinates $(t,\theta)$ on
$(S^1\times[\log Kr_0,\log sH^{-1}])$, then the tension field

\begin{equation}
  |\tau(v)|=r^2|\nabla_e H_e|\leq Cr_0^{-1}
\end{equation}
for $t\in[\log Kr_0,\log sH^{-1}]$. Thus,

\begin{equation}
\int_{S^1\times[t,t+L]}|\tau(v)|^2dtd\theta\leq Cr_0^{-2}
\end{equation}

Let $I_i$ stand for $S^1\times[\log Kr_0+(i-1)L,\log Kr_0+iL]$, and
$N_i$ stand for $I_{i-1}\cup I_i\cup I_{i+1}$. On $\Sigma_n$ we
assume $\log(sH^{-1})-\log(Kr_0)=l_nL$. And like \cite{QT1}, first
we prove that,

\lemma For each $i\in[3,l_n-2]$, there exists a geodesic $\gamma$
such that

\begin{equation}\label{22}
  \int_{I_i}|\widetilde{\nabla}(v-\gamma)|^2dtd\theta\leq
  C(e^{-iL}+e^{-(l_n-i)L})s^2+Cr_0^{-1}
\end{equation}
where $\widetilde{\nabla}$ is the gradient on
$S^1\times[\log(Kr_0),\log(sH^{-1})]$

\proof By Theorem\ref{21}, we have

\begin{equation}\label{35}
  [v]_{C^\alpha(I_i)}\leq \|\widetilde{\nabla}v\|_{L^\infty}\leq C(r_0^{-\frac{1}{2}}+s)
\end{equation}
then if $r_0$ sufficiently large and $s$ sufficiently small, we have
$[v]_{C^\alpha(N_i)}$ is very small.

To apply the Lemma \ref{20} to prove this lemma we choose to points
$P$ and $Q$ on $S^2$(the image of Gauss map) satisfying

\begin{eqnarray}
  &&|P-\frac{1}{2\pi}\int_{(i-1)L\times S^1}vd\theta|\leq C \max_{(i-1)L\times S^1}|v-P|^2\nonumber\\
  &&|Q-\frac{1}{2\pi}\int_{iL\times S^1}vd\theta|\leq C \max_{iL\times S^1}|v-Q|^2
\end{eqnarray}
Note that $S^2$ is compact and smooth, so by (\ref{35}) we can
always find such $P$ and $Q$ and $P$,$Q$ are very close. So there is
a unique geodesic $\gamma_i$ connecting $P$ and $Q$ whose velocity
is sufficiently small.

So if we write down the equation satisfied by $v-\gamma_i$ on
$S^1\times[\log(Kr_0),\log(sH^{-1})]$

\begin{equation}
  \widetilde{\Delta}u+A\cdot\widetilde{\nabla}u+B\cdot u=\tau
\end{equation}
where $u=v-\gamma_i$, we have:

\begin{eqnarray} \label{23}
  &&|A|\leq C(|\widetilde{\nabla} v|+|\widetilde{\nabla}\gamma_i|)\leq
  \delta_0\nonumber\\
  &&|B|\leq C\min\{|\widetilde{\nabla}
  v|^2,|\widetilde{\nabla}\gamma_i|^2\}\leq\delta_0
\end{eqnarray}

If Lemma \ref{20} (C') cannot be used, the only reason is that

\begin{equation}
 \|v-\gamma_i\|_{1,i}\leq C\|\tau\|_{L^2(N_i)}
\end{equation}
which implies

\begin{equation}
  \int_{I_i}|\widetilde{\nabla} (v-\gamma_i)|^2dtd\theta\leq Cr_0^{-2}
\end{equation}
which implied (\ref{22}).

If Lemma \ref{20} (C') can be used, then applying it for
$u=v-\gamma_i$ over $N_i$, we have either

\begin{equation}
  \|u\|_{1,i}<e^{-\frac{1}{2}L}\|u\|_{1,i-1}
\end{equation}
or

\begin{equation}
  \|u\|_{1,i}<e^{-\frac{1}{2}L}\|u\|_{1,i+1}.
\end{equation}
Suppose the first one happens (without loss of generality). Then we
may push this relation to the left because (\ref{23})  hold
regardless of t's position. If the theorem can be used on $N_{j+1}$
but not on $N_j$ for some $j\geq2$, then we have

\begin{equation}
  \|u\|_{1,i}<e^{-\frac{1}{2}(i-j)L}\|u\|_{1,j}\leq
  Ce^{-\frac{1}{2}(i-j)L}r_0^{-1}\leq Cr_0^{-1}.
\end{equation}
If the theorem can be used until $I_2$, then we have

\begin{eqnarray}
  &&e^{\frac{L}{2}}\|u\|_{1,2}\leq\|u\|_{1,1}=(\int_{I_1}u^2dtd\theta)^{\frac{1}{2}}+(\int_{I_1}|\widetilde{\nabla}u|^2dtd\theta)^{\frac{1}{2}}\nonumber\\
  &&\leq(\int_{I_2}u^2dtd\theta)^{\frac{1}{2}}+(\int_{I_1}(u(t,\theta)-u(t+L,\theta))^2dtd\theta)^{\frac{1}{2}}+(\int_{I_1}|\widetilde{\nabla}u|^2dtd\theta)^{\frac{1}{2}}\nonumber\\
\end{eqnarray}
So we have

\begin{eqnarray}
  &&(e^{\frac{L}{2}}-1)\|u\|_{1,2}\leq(\int_{I_1}(\int_0^L|\frac{\partial u}{\partial
  t}(t+s,\theta)|ds)^2dtd\theta)^{\frac{1}{2}}+(\int_{I_1}|\widetilde{\nabla}u|^2dtd\theta)^{\frac{1}{2}}\nonumber\\
  &&\leq \int_0^L(\int_{I_1}|\frac{\partial u}{\partial
  t}(t+s,\theta)|^2dtd\theta)^{\frac{1}{2}}ds+(\int_{I_1}|\widetilde{\nabla}u|^2dtd\theta)^{\frac{1}{2}}\nonumber\\
  &&\leq C (\int_{I_1\cup
  I_2}|\widetilde{\nabla}u|^2dtd\theta)^{\frac{1}{2}}\nonumber\\
  &&\leq C (\int_{I_1\cup
  I_2}|\widetilde{\nabla}v|^2dtd\theta)^{\frac{1}{2}}+C (\int_{I_1\cup
  I_2}|\widetilde{\nabla}\gamma_i|^2dtd\theta)^{\frac{1}{2}}\nonumber\\
  &&\leq C_(r_0^{-\frac{1}{2}}+s)
\end{eqnarray}
So we have the estimate

\begin{equation}
  \|u\|_{1,i}\leq Ce^{-\frac{i-2}{2}L}\|u\|_{1,2}\leq Ce^{-\frac{i}{2}L}(r_0^{-\frac{1}{2}}+s)
\end{equation}

If $\|u\|_{1,i}<e^{-\frac{1}{2}L}\|u\|_{1,i+1}$ happens, we will
have similarly

\begin{equation}
  \|u\|_{1,i}\leq Ce^{-\frac{l_n-i}{2}L}(r_0^{-\frac{1}{2}}+s)
\end{equation}

Finally we get

\begin{equation}
  \|u\|_{1,i}\leq C(e^{-\frac{i}{2}L}+e^{-\frac{l_n-i}{2}L})s+Cr_0^{-\frac{1}{2}}
\end{equation}
which implies (\ref{22}).

Then to get the energy decay, we use the Hopf differential

\begin{equation}
  \Phi=|\partial_tv|^2-|\partial_\theta v|^2-2\sqrt{-1}\partial_tv\cdot\partial_\theta v
\end{equation}

We know that the $L^1$ norm of $\Phi$ is invariant under conformal
change of the coordinates. $(t,\theta)$ is the coordinate of
$A_{Kr_0e^{(i-2)L},Kr_0e^{(i+1)L}}$, we find another coordinate for
it: set $r_i=Kr_0e^{iL}$, then
$(r_i^{-1}A_{Kr_0e^{(i-2)L},Kr_0e^{(i+1)L}}, r_i^{-2}g_e)$ can be
represented as $(A^0_{e^{-2L},e^L},\overline{g})$ , where
$\|\overline{g}-|dx|^2\|_{C^1(A^0_{e^{-2L},e^L})}\leq\varepsilon$.
Assume this Euclidean coordinate is $(x,y)$, so:

\begin{eqnarray}
  \int_{S^1\times[\log Kr_0+(i-1)L,\log
  Kr_0+iL]}|\Phi|dtd\theta=\int_{A^0_{e^{-L},1}}|\Phi|dxdy
\end{eqnarray}

To estimate the right hand side, we use the Cauchy integral formula
on $\Omega=A^0_{e^{-2L},e^L}$, and set
$\Omega^{\prime}=A^0_{e^{-L},1}$, for any $z\in \Omega^{\prime}$

\begin{equation}
  \Phi(v)(z)=\frac{1}{2\pi\sqrt{-1}}\int_{\partial\Omega}\frac{\Phi(w)}{w-z}dw+\frac{1}{2\pi\sqrt{-1}}\int_{\Omega}\frac{\partial\Phi(w)}{\partial\overline{w}}\frac{dw\wedge d\overline{w}}{w-z}
\end{equation}
We know

\begin{equation}
  |\partial_xv|,|\partial_yv|\leq CKr_0e^{iL}|A|\leq
  CKr_0e^{iL}(|x|^{-1}r_0^{-\frac{1}{2}}+r_1^{-1})\leq C(r_0^{-\frac{1}{2}}+se^{-(l_n-i)L})
\end{equation}
so we have:

\begin{equation}
  |\frac{1}{2\pi\sqrt{-1}}\int_{\partial\Omega}\frac{\Phi(w)}{w-z}dw|\leq
  C(r_0^{-1}+s^2e^{-2(l_n-i)L})
\end{equation}

For the second term, notice that by easy calculation

\begin{equation}
  \frac{\partial\Phi(w)}{\partial\overline{w}}=\partial v\cdot
  \overline{\tau}(v)
\end{equation}
where $ \overline{\tau}(v)$ is the tension field under this
coordinate. And

\begin{equation}
 |\overline{\tau}(v)|\leq (Kr_0e^{il})^2|\nabla_eH_e|\leq Cr_0^{-1}
\end{equation}
so we have:

\begin{equation}
  \frac{1}{2\pi\sqrt{-1}}\int_{\Omega}\frac{\partial\Phi(w)}{\partial\overline{w}}\frac{dw\wedge
  d\overline{w}}{w-z}\leq Cr_0^{-1}
\end{equation}

Then we get:

\begin{equation}
  \int_{\Omega^{\prime}}|\Phi|\leq C(r_0^{-1}+s^2e^{-2(l_n-i)L})
\end{equation}

By direct calculation

\begin{eqnarray}
  &&\int_{S^1\times[Kr_0e^{(i-1)L},Kr_0e^{iL}]}|\partial_t
  v|^2dtd\theta\nonumber\\
  &&\leq\int_{S^1\times[Kr_0e^{(i-1)L},Kr_0e^{iL}]}|\Phi|dtd\theta+\int_{S^1\times[Kr_0e^{(i-1)L},Kr_0e^{iL}]}|\partial_{\theta}v|^2dtd\theta\nonumber\\
\end{eqnarray}
and we can get the estimate of
$\int_{S^1\times[Kr_0e^{(i-1)L},Kr_0e^{iL}]}|\partial_{\theta}v|^2dtd\theta$
directly by (\ref{22}). So we get the estimate:

\begin{equation}
  \int_{S^1\times[Kr_0e^{(i-1)L},Kr_0e^{iL}]}|\widetilde{\nabla}v|^2dtd\theta\leq C(e^{-iL}+e^{-(l_n-i)L})s^2+Cr_0^{-1}
\end{equation}

\proposition \label{25}Suppose that$\{\Sigma_n\}$ is a sequence of
constant mean curvature surfaces in a given asymptotically flat end
$(R^3\setminus B_1(0),g)$ and that

\begin{equation}
  \lim_{i\rightarrow \infty}r_0(\Sigma_n)=\infty
\end{equation}

And suppose that

\begin{equation}
\lim_{n\rightarrow \infty} r_0(\Sigma_n)H(\Sigma_n)=0
\end{equation}

Then there exist a large number $K$, a small number $s$ and $n_0$
such that,when $n\geq n_0$,

\begin{equation}
  \max_{I_i}|\widetilde{\nabla}v|\leq C(e^{-\frac{i}{2}L}+e^{-\frac{(l_n-i)}{2}L})s+Cr_0^{-\frac{1}{2}}
\end{equation}
where

\begin{equation}
  I_i=S^1\times[\log(Kr_0(\Sigma_n))+(i-1)L,\log(Kr_0(\Sigma_n))+iL]
\end{equation}
and

\begin{eqnarray}
  i\in[0,l_n] && \log(Kr_0(\Sigma_n))+l_nL=\log(sH^{-1}(\Sigma_n))
\end{eqnarray}

Proof. We just use the interior estimate of the elliptic equation

\begin{equation}
  \widetilde{\Delta}v+|\widetilde{\nabla}v|^2v=\tau
\end{equation}

We know $\|\widetilde{\nabla}v\|_{\infty}\leq
C(r_0^{-\frac{1}{2}}+s)$ , and $\|\tau\|_{\infty}\leq Cr_0^{-1}$.
Assume that :

\begin{equation}
  I_i\subset\subset \widetilde{I}_i
  \subset\subset N_i
\end{equation}
then for some $p>2$

\begin{eqnarray}
  &&\sup_{I_i}|\widetilde{\nabla}v|\leq
  C\|\widetilde{\nabla}v\|_{W^{1,p}(I_i)}\leq
  C(\|v\|_{L^p(\widetilde{I}_i)}+r_0^{-1})\leq
  C(\|v\|_{L^2(N_i)}+r_0^{-1})\nonumber\\
 &&\leq C(e^{-\frac{i}{2}L}+e^{-\frac{(l_n-i)}{2}L})s+Cr_0^{-\frac{1}{2}}
\end{eqnarray}

This analysis improves our understanding of the blowdowns that we
discussed in the previous section. Namely,

\corollary \label{27}Assume the same condition as the above
proposition and in addition $ \lim_{r_0\rightarrow
\infty}\frac{log(r_1)}{r_{0}^{1/4}}=0$. Then the limit plane in
Lemma\ref{17} and Lemma\ref{19} are all orthogonal to the same
vector $a$. In fact, we may choose $s$ small and $i$ large enough so
that,

\begin{equation}
  |v(x)+a|\leq \varepsilon
\end{equation}
for all $x\in \Sigma_n$ and $|x|\leq sH^{-1}(\Sigma_n)$

Proof. We want to prove that

\begin{equation}
  Osc_{B_{sH^{-1}}\cap \Sigma_n}v
\end{equation}
is sufficiently small if $r_0(\Sigma_n)$ large and $s$ small. We
already know that

\begin{equation}
  Osc_{B_{Kr_0}\cap \Sigma_n}v
\end{equation}
is very small from Lemma \ref{17}, so we need only to prove that

\begin{equation}
  Osc_{(B_{sH^{-1}}\setminus B_{Kr_0})\cap \Sigma_n}v
\end{equation}
is small.

From the proposition above we find that

\begin{eqnarray}
  Osc_{(B_{sH^{-1}}\setminus B_{Kr_0})\cap \Sigma_n}v\leq
  \sum_{i=1}^{l_n}Osc_{I_i}v\leq
  C\sum_{i=1}^{l_n}\sup_{I_i}|\widetilde{\nabla}v|\nonumber\\
  \leq
  C\sum_{i=1}^{l_n}((e^{-\frac{i}{2}L}+e^{-\frac{(l_n-i)}{2}L})s+r_0^{-\frac{1}{2}})\leq
  Cs+l_nr_0^{-\frac{1}{2}}
\end{eqnarray}

From $C^{-1}r_1\leq H^{-1}\leq Cr_1$ and the condition
$\lim_{r_0\rightarrow \infty}\frac{\log(r_1)}{r_{0}^{1/4}}=0$, we
have
\begin{equation}
  l_nr_0^{-\frac{1}{2}}=L^{-1}(\log(sH^{-1})-\log(Kr_0))r_0^{-\frac{1}{2}}\leq
  C\frac{\log r_1}{r_0^{\frac{1}{2}}}\rightarrow0
\end{equation}
as $r_0\rightarrow\infty$, so we prove the lemma.

\corollary \label{28}Assume the same condition as Proposition
\ref{25}. Let $v_n=v(p_n)$ for some $p_n\in I_{\frac{l_n}{2}}$. Then

\begin{equation}
  \sup_{I_i}|v-v_n|\leq
  C(e^{-\frac{1}{2}iL}+e^{-\frac{1}{4}l_nL})s+l_nr_0^{-\frac{1}{2}}
\end{equation}
for $i\in [0,\frac{1}{2}l_n]$

\begin{equation}
  \sup_{I_i}|v-v_n|\leq
  C(e^{-\frac{1}{4}l_nL}+e^{-\frac{1}{2}(l_n-i)L})s+l_nr_0^{-\frac{1}{2}}
\end{equation}
for $i\in [\frac{1}{2}l_n,l_n]$

\section{Harmonic Coordinates}
We assume that the metric $g$ can be expanded in the coordinate
$\{x_i\}$ as
\begin{eqnarray}
  g_{ij}=\delta_{ij}+h_{ij}\nonumber
  =\delta_{ij}+h^1_{ij}(\theta)/r+Q
\end{eqnarray}
where $\theta$ is the coordinate on the unit sphere $S^2$, and
$h^1_{ij}(\theta)$ is a function extended constantly along the
radius direction. And $Q$ satisfies

\begin{equation}
  \sup r^{2+k}|\partial^kQ|\leq C
\end{equation}
for $k=0,1,\cdots,5$

First, note that:
\begin{eqnarray}
  &&\Delta_gx_k=\frac{1}{\sqrt{g}}\frac{\partial}{\partial x_i}(\sqrt{g}g^{ij}\frac{\partial}{\partial
  x_j}x_k)\nonumber\\&&=\frac{\partial}{\partial
  x_i}g^{ik}+\frac{1}{2}g^{ik}g^{mn}g_{mn,i}\nonumber\\
  &&=-g^{mn}\Gamma_{mn}^k=O(|x|^{-2})
\end{eqnarray}

Now our aim is to find asymptotically harmonic coordinate, i.e. some
coordinate ${y^i}$ such that $\Delta_gy^k=O(|x|^{-3})$

\begin{eqnarray}
  &&\Delta_gx^k=-g^{jl}g^{ik}\frac{1}{2}(\frac{\partial}{\partial x_j}h_{li}+\frac{\partial}{\partial x_l}h_{ji}-\frac{\partial}{\partial
  x_i}h_{jl})\nonumber\\
  &&=-g^{jl}g^{ik}\frac{1}{2}(
  r^{-2}((h^1_{li,j}(\theta)-h^1_{li}(\theta)\frac{x_j}{r})\nonumber\\
  &&+(h^1_{ji,l}(\theta)-h^1_{ji}(\theta)\frac{x_l}{r})-(h^1_{jl,i}(\theta)-h^1_{jl}(\theta)\frac{x_i}{r})))+\partial Q\nonumber\\
  &&=-g^{jl}g^{ik}\frac{1}{2}r^{-2}f^1_{lij}(\theta)+O(|x|^{-3})
\end{eqnarray}

We also know that $g^{ij}=\delta^{ij}-h^1_{ij}(\theta)/r+O(r^{-2})$

Then :
\begin{equation}
   \Delta_gx^k=-\frac{1}{2}r^{-2}f^1_{jkj}(\theta)+O(r^{-3})
\end{equation}

Suppose $0=\xi_0>\xi_1\geq\xi_2\geq\cdots$ are the eigenvalues of
$\Delta|_{S^2}$, and $A_n(\theta)$ are the corresponding orthonormal
eigenvectors.

Set:
\begin{equation}
  y^k=x^k+\sum_{n=0}^\infty f^k_n(r)A_n(\theta)
\end{equation}

We have:
\begin{equation}
  \Delta_g y^k=\Delta_g x^k+\sum_{n=0}^\infty
  \Delta_{R^3}(f_n^k(r)A_n(\theta))+\sum_{n=0}^\infty(\Delta_g-\Delta_{R^3})(f_n^k(r)A_n(\theta))
\end{equation}

Solve the equation:
\begin{eqnarray}
\Delta_g x^k+\sum_{n=0}^\infty
  \Delta_{R^3}(f^k_n(r)A_n(\theta))=O(|x|^{-3})
\end{eqnarray}

 Assume
\begin{equation}\label{9}
  \frac{1}{2}f_{jkj}^1(\theta)=\sum_{n=0}^\infty \lambda^k_nA_n(\theta)
\end{equation}

so we have:
\begin{eqnarray}
\sum_{n=0}^\infty
\Delta_{R^3}(f_n^k(r)A_n(\theta))=r^{-2}\sum_{n=0}^\infty
\lambda^k_nA_n(\theta)\\
\frac{1}{r^2}(2rf_n^{k\prime}+r^2f_n^{k\prime\prime}+f_n^k(r)\xi_n)=\lambda_n^k
, n=0,\cdots ,\infty
\end{eqnarray}

\begin{eqnarray}
&n=0 , &f_0^k=\lambda_0^k\log(r)\\
&n>0 , &f_n^k=\frac{\lambda^k_n}{\xi_n}
\end{eqnarray}

and this solution satisfies that:
\begin{eqnarray}
\sum_{n=0}^\infty(\Delta_g-\Delta_{R^3})(f_n^k(r)A_n(\theta))=O(|x|^{-3})
\end{eqnarray}

 so if
\begin{equation}\label{8}
y^k=x^k+\frac{1}{2\sqrt{\pi}}\lambda_0^k\log r+\sum_{n=1}^\infty
\frac{\lambda^k_n}{\xi_n} A_n(\theta)
\end{equation}

then we must have:
\begin{equation}
    \Delta y^k=O(|x|^{-3})
\end{equation}

 Note that
\begin{eqnarray}
\Delta|_{S^2}\sum_{n=1}^\infty \frac{\lambda^k_n}{\xi_n}
A_n(\theta)=\sum_{n=1}^\infty\lambda^k_nA_n(\theta)=\frac{1}{2}f_{jkj}^1(\theta)-\frac{1}{2}\overline{f_{jkj}^1(\theta)}
\end{eqnarray}
where $\overline{f_{jkj}^1(\theta)}$ is its mean value on the unit
sphere.

Set
\begin{equation}
    g_k^1(\theta)=\sum_{n=1}^\infty \frac{\lambda^k_n}{\xi_n}
A_n(\theta)=\Delta^{-1}(\frac{1}{2}f^1_{jkj}(\theta)-\frac{1}{2}\overline{f^1_{jkj}(\theta)})
\end{equation}

\begin{equation}\label{7}
\frac{\partial y^k}{\partial
x^i}=\delta_{ik}+\frac{\lambda^k_0}{2\sqrt{\pi}}\frac{1}{r}\frac{x^i}{r}+g_k^1(\theta)_i\frac{1}{r}
\end{equation}
\begin{equation}
\frac{\partial x^i}{\partial y^k}=\delta_{ik}+O(|x|^{-1})
\end{equation}

So we get:
\begin{equation}
\widetilde{g}_{ij}=g(\frac{\partial}{\partial
y^i},\frac{\partial}{\partial y^j})=\delta_{ij}+O(|x|^{-1})
\end{equation}

Suppose
\begin{equation}
\widetilde{g}_{ij}=\delta_{ij}+\widetilde{h}_{ij}
\end{equation}

Now I want to discuss the ellipticity of $\widetilde{h}_{ij}$
\begin{equation}
\widetilde{h}_{ij}=h_{ij}-\frac{1}{2r\sqrt{\pi}}(\lambda_0^i\frac{x^j}{r}+\lambda_0^j\frac{x^i}{r})-\frac{(g^1_{i,j}(\theta)+g^1_{j,i}(\theta))}{r}
\end{equation}
Where $g^1_{i,j}(\theta)$ denotes the constant extention along the
radius direction of function$\frac{\partial
g^1_{i}(\theta)}{\partial x_j}|_{S^2}$

\example: For the metric
 $g_{ij}=\delta_{ij}+\frac{\delta_{ij}}{r}$, we have:
 \begin{equation}
   \Delta_{g}x^k=-\frac{1}{2}\frac{x^k}{r^3}+O(|x|^{-3})
 \end{equation}
We know that on $S^2$, we have $\Delta|_{S^2}x^k=-2x^k$ . So if we
let:
\begin{equation}
  y^k=x^k-\frac{1}{4}\frac{x^k}{r}
\end{equation}

We have $\Delta_gy^k=O(|x|^{-3})$ , then:
\begin{equation}
  \frac{\partial y^k}{\partial x^i}=\delta_{ki}-\frac{1}{4}(\frac{\delta_{ki}}{r}-\frac{x^kx^i}{r^3})
\end{equation}
\begin{equation}
  \widetilde{h}_{ij}=\frac{3\delta_{ij}}{2r}-\frac{x^ix^j}{2r^3}+O(r^{-2})
\end{equation}

\lemma Suppose in some coordinate $\{x^i\}$ ,
$g_{ij}=\delta_{ij}+h^1_{ij}(\theta)/r+Q$ , then for any $m>2$ there
exists $\varepsilon>0$ , if
  $\|h^1_{ij}(\theta)-\delta_{ij}(\theta)\|_{W^{m,2}(S^2)}\leq
\varepsilon$  then in the asymptotically harmonic coordinate
$\{y^i\}$we get above , we have
\begin{equation}
  \widetilde{g}_{ij}=\delta_{ij}+\widetilde{h}_{ij}
\end{equation}
where $\widetilde{h}_{ij}=O(|y|^{-1})$ , and $|y|\widetilde{h}_{ij}$
is uniformly elliptic.

Proof: We know easily from (\ref{7}) that
$\widetilde{h}_{ij}=O(|x|^{-1})$ and that
$\lim_{|x|\rightarrow\infty}\frac{|y|}{|x|}=1$ ,then
$\widetilde{h}_{ij}=O(|y|^{-1})$ . So we need only to prove that
$|y|\widetilde{h}_{ij}$ is uniformly elliptic.

First we know from
$\|h^1_{ij}(\theta)-\delta_{ij}(\theta)\|_{W^{m,2}(S^2)}\leq
\varepsilon$  that
\begin{eqnarray}
  \|\frac{1}{2}f^1_{jkj}(\theta)-\frac{1}{2}\frac{x^k}{r}\|_{W^{m-1,2}(S^2)}\leq
C \varepsilon
\end{eqnarray}

Note that $\frac{1}{2}f_{jkj}^1(\theta)=\sum_{n=0}^\infty
\lambda^k_nA_n(\theta)$ and $x^k$ is an eigenvector of
$\Delta_{S^2}$ , so we can assume that $A_1(\theta)=C_kx^k|_{S^2}$
without loss of generality.
\begin{eqnarray}
  \|\lambda^k_0A_0(\theta)+(\lambda_1^kC_k-\frac{1}{2})x^k+\sum_{n=2}^{\infty}\lambda^k_nA_n(\theta)\|_{W^{m-1,2}(S^2)}\leq \varepsilon
\end{eqnarray}
so we get
\begin{eqnarray}
  |\lambda^k_0|\leq\varepsilon , (\lambda_1^kC_k-\frac{1}{2})\leq\varepsilon , \sum_{n=2}^\infty (|\xi_n|^{\frac{m-1}{2}}\lambda^k_n)^2\leq\varepsilon
\end{eqnarray}

Note that from (\ref{8})
\begin{eqnarray}
  \frac{\partial y^k}{\partial
x^i}=\delta_{ik}+\frac{\lambda^k_0}{2\sqrt{\pi}}\frac{1}{r}\frac{x^i}{r}-\frac{1}{2}(\frac{1}{2}\pm\varepsilon)(\frac{\delta_{ik}}{r}-\frac{x_ix_k}{r^3})
+\sum_{n=2}^\infty \frac{\lambda^k_n}{\xi_n}\frac{\partial A_n(\theta)}{\partial x_i}
\end{eqnarray}
where the last term on the right can be estimated, for some $p>0$
\begin{eqnarray}
  &&|\sum_{n=2}^\infty \frac{\lambda^k_n}{\xi_n}\frac{\partial A_n(\theta)}{\partial x_i}|\leq \sum_{n=2}^\infty \frac{|\lambda^k_n|}{|\xi_n|}\frac{|\nabla_{S^2}A_n(\theta)|}{r}\nonumber\\
  &&\leq \sum_{n=2}^\infty \frac{|\lambda^k_n|}{|\xi_n|}\frac{\|A_n(\theta)\|_{W^{2+p,2}}}{r}\nonumber\\
  &&\leq \sum_{n=2}^\infty \frac{|\lambda^k_n|}{|\xi_n|}\frac{|\xi_n|^{1+\frac{p}{2}}\|A_n(\theta)\|_{L^2}}{r}\nonumber\\
  &&\leq \frac{1}{r}\sum_{n=2}^\infty |\lambda^k_n||\xi_n|^{\frac{m-1}{2}}|\xi_n|^{\frac{p-m+1}{2}}\nonumber\\
  &&\leq \frac{1}{r}(\sum_{n=2}^\infty (|\lambda^k_n||\xi_n|^{\frac{m-1}{2}})^2)^\frac{1}{2}(\sum_{n=2}^\infty|\xi_n|^{p-m+1})^{\frac{1}{2}}
\end{eqnarray}
let $p=\frac{m-2}{2}$ , then from $\xi_n=O(n)$ we have
\begin{equation}
  \sum_{n=2}^\infty|\xi_n|^{p-m+1}\leq C
\end{equation}
so
\begin{equation}
 |\sum_{n=2}^\infty \frac{\lambda^k_n}{\xi_n}\frac{\partial A_n(\theta)}{\partial x_i}|\leq \frac{C\varepsilon}{r}
\end{equation}
then we have:
\begin{equation}
   \frac{\partial y^k}{\partial
x^i}=\delta_{ik}-\frac{1}{4}(\frac{\delta_{ik}}{r}-\frac{x_ix_k}{r^3})+\frac{C\varepsilon}{r}
\end{equation}
so we can deduce that:
\begin{equation}
  \widetilde{h}_{ij}=h_{ij}+\frac{\delta_{ij}}{2r}-\frac{x^ix^j}{2r^3}+\frac{C\varepsilon}{r}
\end{equation}
because $|h^1_{ij}(\theta)-\delta_{ij}(\theta)|_{W^{m,2}(S^2)}\leq
\varepsilon$, we have $rh_{ij}$ is uniformly elliptic. And the
eigenvalues of $\frac{x^ix^j}{r^2}$ are between $0$ and $1$, so
$|y|\widetilde{h}_{ij}$ is uniformly elliptic from
$\lim_{r\rightarrow \infty}\frac{|y|}{r}=1$ for $\varepsilon$
sufficiently small.

So all the analysis in Section 2,3,4 can be done in the asymptotically harmonic coordinate $\{y_i\}$.

\lemma
In the asymptotically harmonic coordinate $\{y^i\}$, we have that
\begin{eqnarray}
  -\frac{1}{2}\Delta_{g}\log|\widetilde{g}|=R(g)+O(|y|^{-4})
\end{eqnarray}
Proof. From direct calculation we have
\begin{eqnarray}
  R(g)=\widetilde{g}^{jk}\widetilde{g}^{il}\widetilde{g}_{ml}(\frac{\partial\widetilde{\Gamma}^m_{jk}}{\partial y^i}
  -\frac{\partial\widetilde{\Gamma}^m_{ik}}{\partial y^j})+O(|y|^{-4})
\end{eqnarray}
\begin{eqnarray}
  &&\widetilde{g}^{jk}\widetilde{g}^{il}\widetilde{g}_{ml}\frac{\partial\widetilde{\Gamma}^m_{jk}}{\partial y^i}
  =\widetilde{g}^{il}\widetilde{g}_{ml}\frac{\partial (\widetilde{g}^{jk}\widetilde{\Gamma}^m_{jk})}{\partial y^i}+O(|y|^{-4})\nonumber\\
  &&=-\widetilde{g}^{il}\widetilde{g}_{ml}\frac{\partial\Delta_gy^m}{\partial y^i}+O(|y|^{-4})=O(|y|^{-4})
\end{eqnarray}
\begin{eqnarray}
    &&-\widetilde{g}^{jk}\widetilde{g}^{il}\widetilde{g}_{ml}\frac{\partial\widetilde{\Gamma}^m_{ik}}{\partial y^j}
    =-\frac{1}{2}\widetilde{g}^{jk}\widetilde{g}^{ip}\frac{\partial^2\widetilde{g}_{ip}}{\partial y^j\partial y^k}+O(|y|^{-4})\nonumber\\
    &&= -\frac{1}{2}\Delta_g\log|\widetilde{g}|+O(|y|^{-4})
\end{eqnarray}
so we prove the lemma.

\corollary If in addition $R=O(|x|^{-3-\tau})$ for some $\tau>0$ ,
then in the asymptotically harmonic coordinate $\{y^i\}$ ,
we have
\begin{equation}
  \sum_{i=1}^3\widetilde{h}_{ii}=8m(g)/|y|+o(|y|^{-1-\frac{\tau}{2}})
\end{equation}

Proof: First we know that

\begin{equation}
  \lim_{|x|\rightarrow\infty}\frac{|y|}{|x|}=1,
\end{equation}then from the lemma above that in the coordinate $\{y^i\}$ , we
have
\begin{eqnarray}
  \Delta_{g}\log|\widetilde{g}|=O(|y|^{-3-\tau})
\end{eqnarray}

We know that
\begin{equation}
  \log|\widetilde{g}|=O(|y|^{-1})
\end{equation}

From the theory of harmonic functions in $R^n$ , we have there exist
some constant $C$ such that:
\begin{equation}
  \log|\widetilde{g}|=\frac{C}{|y|}+o(|y|^{-1-\frac{\tau}{2}})
\end{equation}

From Bartnik's result , we know the mass is invariant under the
change of coordinates because $R(g)\in L^1$ .

\begin{eqnarray}
   m(g)=\lim_{R\rightarrow\infty}\frac{1}{16\pi}\int_{s_R}
   (\widetilde{h}_{ij,j}-\widetilde{h}_{jj,i})v_g^id\mu
\end{eqnarray}
Now we have
\begin{eqnarray}
 &&\widetilde{g}_{ik,k}-\frac{1}{2}\widetilde{g}_{kk,i}=\widetilde{g}^{ij}\widetilde{g}^{kl}(\widetilde{g}_{jk,l}-\frac{1}{2}\widetilde{g}_{kl,j})+O(|y|^{-3})\nonumber\\
 &&=-\Delta_g
  y^i+O(|y|^{-3})=O(|y|^{-3})
\end{eqnarray}

So we have:
\begin{eqnarray}
  m(g)&&=\lim_{R\rightarrow\infty}\frac{1}{16\pi}\int_{s_R}
  (-\frac{1}{2}\widetilde{h}_{jj,i})v_g^id\mu\nonumber\\
  &&=-\lim_{R\rightarrow\infty}\frac{1}{32\pi}\int_{s_R}\frac{\partial\log|\widetilde{g}|}{\partial
  y^i}v_g^id\mu\nonumber\\
  &&=\lim_{R\rightarrow\infty}\frac{1}{32\pi}\int_{s_R}\frac{Cy^i}{|y|^3}v_g^id\mu\nonumber\\
  &&=\frac{C}{8}
\end{eqnarray}

So we get the result by easy calculation .

\begin{remark}
  In fact we can replace the constraint equation by the condition
  \begin{equation}
    R=O(|x|^{-3-\tau})
  \end{equation}
  for some $\tau>0$.
\end{remark}
\section{Proof of the Theorem}

Now let's prove Theorem \ref{26}.

First recall that, for any surface $\Sigma$ embedded in
$\mathbb{R}^3$ and any given vector  $b\in \mathbb{R}^3$, one has
\begin{equation}
  \int_\Sigma H_e<v_e\cdot b>_ed\mu_e=0
\end{equation}
where $H_e$ and $v_e$ denote the mean curvature
and normal vector field with respect to the Euclidean metric.

On the other hand , if $\Sigma$ is a constant mean curvature surface
in the asymptotically flat end , then
\begin{equation}
  \int_\Sigma H<v_e\cdot b>_ed\mu_e=0
\end{equation}

So we have
\begin{eqnarray}
  \int_\Sigma (H-H_e)<v_e\cdot b>_e d\mu_e=0
\end{eqnarray}

From now on , our calculation is in the coordinate $\{x^i\}$ ,which
is assumed to be the asymptotically harmonic coordinate. We have
calculated $H-H_e$, so we have
\begin{eqnarray}\label{10}
  &&\int_\Sigma (H-H_e)<v_e\cdot b>_ed\mu_e=\int_\Sigma(-f^{ik}h_{kl}f^{lj}A_{ij}+\frac{1}{2}Hv^iv^jh_{ij}-f^{ij}v^l\overline{\nabla}_ih_{jl}\nonumber\\
  &&+\frac{1}{2}f^{ij}v^l\overline{\nabla}_lh_{ij}\pm
  C|h||\overline{\nabla} h|\pm C|h|^2|A|)<v_e\cdot b>_ed\mu_e
\end{eqnarray}

We assume that there exists a sequence of constant mean curvature
surfaces $\Sigma_n$ with
\begin{eqnarray}
  \lim_{n\rightarrow\infty}r_0(\Sigma_n)=\infty && \lim_{n\rightarrow\infty}H(\Sigma_n)r_0(\Sigma_n)=0
\end{eqnarray}
otherwise we have get the result from the uniqueness theorem of
Lan-Hsuan Huang. So we can choose $s$ sufficiently small and $K$
sufficiently large with $sH^{-1}>Kr_0$ for $r_0$ sufficiently large.

We know that
\begin{equation}
  |h|=O(|x|^{-1}) , |\overline{\nabla} h|=O(|x|^{-2}) , |A|\leq CH+C|{\AA}|
\end{equation}
from the estimate
\begin{equation}
  |{\AA}|\leq r_0^{-\frac{1}{2}}O(|x|^{-1})
\end{equation}
we have
\begin{eqnarray}
  &&|\int_\Sigma(\pm C|h||\overline{\nabla} h|\pm C|h|^2|A|)<v_e\cdot b>_ed\mu_e|\leq C\int_\Sigma
  (H|x|^{-2}+|x|^{-3})\nonumber\\
&&=O(r_0^{-1})
\end{eqnarray}
by the estimates in Section 2.





Now we calculate other terms in (\ref{10})
\begin{eqnarray}
  &&\int_{\Sigma_n} -f^{ij}v^l(\overline{\nabla}_ih_{jl})v^mb^md\mu_e\nonumber\\
  &=&\frac{1}{2}\int_{\Sigma_n}(f^{ij}h_{jk}f^{kl}A_{li}-Hv^jv^lh_{jl})v^mb^md\mu_e+\frac{1}{2}\int_{\Sigma_n}
  f^{ij}v^lh_{jl}A_{ik}f^{km}b^md\mu_e\nonumber\\
  &&-\frac{1}{2}\int_{\Sigma_n}
  f^{ij}v^l(\overline{\nabla}_ih_{jl})v^mb^md\mu_e
\end{eqnarray}
because $d\mu_e=(1+O(r^{-1}))d\mu$ , $v_e=(1+O(r^{-1}))v$ and
$<v_e\cdot b>_e=<v\cdot b>_g+O(r^{-1})$.

So we have
\begin{eqnarray}\label{11}
  \int_{\Sigma_n} (H-H_e)<v_e\cdot b>_ed\mu_e=\int_{\Sigma_n}-\frac{1}{2}f^{ik}h_{kl}f^{lj}A_{ij}v^mb^m+
  f^{ij}v^lh_{jl}A_{ik}f^{km}b^m\nonumber\\
  -\frac{1}{2}f^{ij}v^l\overline{\nabla}_ih_{jl}v^mb^m
  +\frac{1}{2}f^{ij}v^l\overline{\nabla}_lh_{ij}v^mb^m+O(r_0^{-1})d\overline{\mu}\nonumber\\
\end{eqnarray}
Note that
\begin{eqnarray}
  A_{ij}={\AA}_{ij}+\frac{f_{ij}}{2}H           , \ \ \            \sup|{\AA}|\leq r_0^{-\frac{1}{2}}O(|x|^{-1})
\end{eqnarray}

So we have
\begin{eqnarray}
  &&\int_{\Sigma_n} (H-H_e)<v_e\cdot b>_ed\mu_e=\int_{\Sigma_n}
  -\frac{H}{4}f^{kl}h_{kl}v^mb^m+\frac{H}{4}f^{jm}h_{jl}v^lb^m\nonumber\\
  &&+\frac{1}{2}f^{ij}(\overline{\nabla}_lh_{ij})v^lv^mb^m-\frac{1}{2}f^{ij}(\overline{\nabla}_ih_{jl})v^lv^mb^m\nonumber\\
  &&\pm C\int_{\Sigma_n}
  |x|^{-2}r_0^{-\frac{1}{2}}+O(r_0^{-1})
\end{eqnarray}


In this case we calculate
\begin{eqnarray}
\int_{\Sigma_n}|x|^{-2}r_0^{-\frac{1}{2}}d\mu_e
\end{eqnarray}
We divide the integral into three parts:
\begin{eqnarray}
  \int_{\Sigma_n}|x|^{-2}r_0^{-\frac{1}{2}}=\int_{\Sigma_n\cap
  B^c_{sH^{-1}}(0)}+\int_{\Sigma_n\cap B_{Kr_0}(0)}+\int_{\Sigma_n\cap (B_{sH^{-1}}\setminus
  B_{Kr_0})}|x|^{-2}r_0^{-\frac{1}{2}}.
\end{eqnarray}

Then by the blowdown results in Section 3 we have
\begin{eqnarray}
  \int_{\Sigma_n\cap
  B^c_{sH^{-1}}(0)}|x|^{-2}r_0^{-\frac{1}{2}}d\mu_e=\int_{\widetilde{\Sigma}_n\cap
  B^c_s(0)}|\widetilde{x}|^{-2}r_0^{-\frac{1}{2}}d\widetilde{\mu}\leq Cr_0^{-\frac{1}{2}}
\end{eqnarray}

\begin{eqnarray}
  \int_{\Sigma_n\cap
  B_{Kr_0}(0)}|x|^{-2}r_0^{-\frac{1}{2}}d\mu_e=\int_{\widehat{\Sigma}_n\cap
  B_k(0)}|\widehat{x}|^{-2}r_0^{-\frac{1}{2}}d\widehat{\mu}\leq
  Cr_0^{-\frac{1}{2}}
\end{eqnarray}

\begin{eqnarray}
  &&\int_{\Sigma_n\cap (B_{sH^{-1}}\setminus
  B_{Kr_0})}|x|^{-2}r_0^{-\frac{1}{2}}d\mu_e=\sum_{i=0}^n\int_{\Sigma_n\cap (B_{Kr_0e^{4iL}}\setminus
  B_{Kr_0e^{4(i-1)L}})}|x|^{-2}r_0^{-\frac{1}{2}}d\mu_e\nonumber\\
  &&\leq C\sum_{i=0}^n\int_{B_{e^{4L}}\setminus
  B_1}|\overline{x}|^{-2}r_0^{-\frac{1}{2}}d\overline{\mu}
  \leq Cr_0^{-\frac{1}{2}}l_nL
\end{eqnarray}
where $e^{l_nL}Kr_0=sH^{-1}$

so if
\begin{equation}
  \lim_{r_0\rightarrow0}\frac{|\log H|}{r_0^{\frac{1}{2}}}=0
\end{equation}
in other words
\begin{equation}
  \lim_{r_0\rightarrow0}\frac{|\log r_1|}{r_0^{\frac{1}{2}}}=0
\end{equation}
 we have
\begin{equation}
  \int_{\Sigma}|x|^{-2}r_0^{-\frac{1}{2}}d\overline{\mu}\rightarrow0
\end{equation}
as $r_0\rightarrow \infty$

From the property of the asymptotically harmonic coordinate
\begin{eqnarray}
 g^{ij}h_{ij}=\frac{8m(g)}{r}+o(r^{-1-\frac{\tau}{2}})\\
 g^{kl}(g_{ik,l}-\frac{1}{2}g_{kl,i})=O(|x|^{-3})
\end{eqnarray}
\begin{eqnarray}
&&\int_{\Sigma_n}
-\frac{H}{4}f^{kl}h_{kl}v^mb^m+\frac{H}{4}f^{jm}h_{jl}v^lb^m
  +\frac{1}{2}f^{ij}(\overline{\nabla}_lh_{ij}-\overline{\nabla}_ih_{jl})v^lv^mb^m\nonumber\\
  &&=\int_{\Sigma_n}-\frac{H}{4}g^{kl}h_{kl}v^mb^m+\frac{H}{4}g^{jm}h_{jl}v^lb^m\nonumber\\
  &&+\frac{1}{2}g^{ij}(\overline{\nabla}_lh_{ij}-\overline{\nabla}_ih_{jl})v^lv^mb^m
  +O(|r_0|^{-1})\nonumber\\
  &&=-2m(g)\int_{\Sigma_n} (\frac{H}{r}<v_e\cdot b_e>_e+\frac{<x\cdot v_e>_e<v_e\cdot
  b>_e}{r^3})\nonumber\\
  &&+\int_{\Sigma_n} \frac{H}{4}h_{ml}v^lb^m+o(1).\nonumber\\
\end{eqnarray}

So we have:
\begin{eqnarray}
  \lim_{n\rightarrow\infty}(-2m(g)\int_{\Sigma_n} (\frac{H}{r}<v_e\cdot b>_e+\frac{<x\cdot v_e>_e<v_e\cdot
  b>_e}{r^3})+\int_{\Sigma_n} \frac{H}{4}h_{ml}v^lb^m)=0\nonumber\\
\end{eqnarray}

Note that:
\begin{equation}
  h_{ml}v^l=(h_{ml}-\frac{tr(h)}{2}\delta_{ml})v^l+\frac{tr(h)}{2}v^m
\end{equation}
where   $tr(h)=g^{ij}h_{ij}$

Assume that the three eigenvalues of $h_{ml}$ are
\begin{equation}
  \lambda_1\geq \lambda_2\geq \lambda_3\geq 0
\end{equation}

For $p\in \Sigma$ fixed , choose coordinate properly such that
\begin{equation}
 h_{ml}-\frac{tr(h)}{2}\delta_{ml}
\end{equation}
can be written as

\begin{eqnarray}
  \left(
  \begin{array}{ccc}
    \lambda_1-\frac{tr(h)}{2} & 0 & 0 \\
    0 & \lambda_2-\frac{tr(h)}{2} & 0 \\
    0 & 0 & \lambda_3-\frac{tr(h)}{2} \\
  \end{array}
\right)
\end{eqnarray}

Assume $v=(\widetilde{v}^1,\widetilde{v}^2,\widetilde{v}^3)$,and
$(\widetilde{v}^1)^2+(\widetilde{v}^2)^2+(\widetilde{v}^3)^2=1$ .
Then we have
\begin{eqnarray}
  \sum_{i=1}^3((\lambda_i-\frac{tr(h)}{2})\widetilde{v}^i)^2=\frac{(tr(h))^2}{4}-\sum_{i=1}^3\lambda_i(tr(h)-\lambda_i)(\widetilde{v}^i)^2
\end{eqnarray}

Because of the uniformly ellipticity we have there exists $C>0$ ,
such that
\begin{equation}
  \frac{th(h)}{C}\leq \lambda_3\leq \lambda_2 \leq \lambda_1\leq
  (1-\frac{1}{C})th(h)
\end{equation}
so
\begin{equation}
\lambda_i(th(h)-\lambda_i)\geq \frac{1}{C}(1-\frac{1}{C})(tr(h))^2
\end{equation}

hence
\begin{equation}
  \sum_{i=1}^3((\lambda_i-\frac{tr(h)}{2})\widetilde{v}^i)^2\leq
  (\frac{1}{4}-\frac{1}{C}(1-\frac{1}{C}))(tr(h))^2
\end{equation}
\begin{eqnarray}
  &&\int_{\Sigma_n} \frac{H}{4}h_{ml}v^lb^m=\int_{\Sigma_n} \frac{H}{4}(\frac{tr(h)}{2}<v_e\cdot
  b>_e+(h_{ml}-\frac{tr(h)}{2}\delta_{ml})v^lb^m)\nonumber\\
  &&\leq \int_{\Sigma_n} \frac{H tr(h)}{4}(\frac{1}{2}<v_e\cdot
  b>_e+\sqrt{\frac{1}{4}-\frac{1}{C}(1-\frac{1}{C})})\nonumber\\
  &&=\int_{\Sigma_n} \frac{H m(g)}{r}(<v_e\cdot
  b>_e+1-\frac{2}{C})
\end{eqnarray}
so we have
\begin{eqnarray}
  &&\int_{\Sigma_n}(H-H_e)<v_e\cdot b>_e\leq -m\int_{\Sigma_n} \frac{H}{r}<v_e\cdot
  b>_e\nonumber\\&&+\frac{2}{r^3}<x\cdot v_e>_e<v_e\cdot
  b>_ed\mu_e
  +(1-\frac{2}{C})m(g)\int_{\Sigma_n}\frac{H}{r}d\mu_e+o(1)\nonumber\\
\end{eqnarray}
as $n\rightarrow\infty$

From Lemma \ref{18}, we have $\frac{H}{2}\Sigma_n$ subconverges to
some sphere $S^2_1(a)$ with $|a|=1$. Now we choose $b=-a$. Then from
the calculation in \cite{QT}, we have
\begin{eqnarray}
  -m(g)\int_{\Sigma_n} \frac{H}{r}<v_e\cdot b>_e&\rightarrow& -\frac{8}{3}\pi m(g)\\
  -m(g)\int_{\Sigma_n} \frac{2}{r^3}<x\cdot v_e>_e<v_e\cdot b>_e&\rightarrow& -\frac{16}{3}\pi m(g)\label{30}\\
  (1-\frac{2}{C})m(g)\int_{\Sigma_n} \frac{H}{r}&\rightarrow& (1-\frac{2}{C})8\pi m(g)
\end{eqnarray}
as  $n\rightarrow\infty$

Because there is a little difference from \cite{QT},we prove them
again. We notice from Lemma \ref{18}, we have $\frac{H}{2}\Sigma_n$
subconverges to some sphere $S_1(a)$ with $|a|=1$, and the first and
third integral converges to $-m(g)\int_{S_1(a)} \frac{2}{r}<v_e\cdot
b>_e= -\frac{8}{3}\pi m(g)$ and $(1-\frac{2}{C})m(g)\int_
{S_1(a)}\frac{2}{r}= (1-\frac{2}{C})8\pi m(g)$ respectively.

To deal with the (\ref{30}), first we notice that

\begin{equation}
  \int_{S^2(a)}\frac{2}{r^3}<x\cdot v_e>_e<v_e\cdot
  b>_ed\mu_e=\frac{4}{3}\pi
\end{equation}  then we break up the integral (\ref{30}) into three parts.
\begin{eqnarray}
  &&\int_{\Sigma_n} \frac{2}{r^3}<x\cdot v_e>_e<v_e\cdot b>_ed\mu_e\nonumber\\
  &&=\int_{\Sigma_n\cap
  B^c_{sH^{-1}}(0)}+\int_{\Sigma_n\cap B_{Kr_0}(0)}+\int_{\Sigma_n\cap B_{sH^{-1}}\setminus B_{Kr_0}}\frac{2}{r^3}<x\cdot v_e>_e<v_e\cdot
  b>_ed\mu_e\nonumber\\
\end{eqnarray}

Then

\begin{eqnarray}
 &&\lim_{n\rightarrow\infty}\int_{\Sigma_n\cap B^c_{sH^{-1}}(0)}\frac{2}{r^3}<x\cdot v_e>_e<v_e\cdot
  b>_ed\mu_e\nonumber\\
  &&=\int_{S^2(a)\cap B^c_s}\frac{2}{r^3}<x\cdot v_e>_e<v_e\cdot
  b>_ed\mu_e
\end{eqnarray}
and

\begin{eqnarray}
  &&\lim_{n\rightarrow\infty}\int_{\Sigma_n\cap B_{Kr_0}(0)}\frac{2}{r^3}<x\cdot v_e>_e<v_e\cdot
  b>_ed\mu_e\nonumber\\
  &&=\int_{P\cap B_K(0)}\frac{2}{r^3}<x\cdot v_e>_e<v_e\cdot
  b>_ed\mu_e,
\end{eqnarray}
where $P$ is the limit plane in Lemma \ref{17}. From
Corollary\ref{27}, we know the normal vector of $P$ is $v_e$. Then
due to an easy calculation we know

\begin{eqnarray}
  \int_P\frac{2}{r^3}<x\cdot v_e>_e<v_e\cdot
  b>_ed\mu_e=4\pi
\end{eqnarray}

From the divergence theorem we have

\begin{equation}
  \int_{\Sigma_n}\frac{2}{r^3}<x\cdot v_e>_ed\mu_e=8\pi
\end{equation}
for any $n$ and

\begin{equation}
  \int_{S^2(a)}\frac{2}{r^3}<x\cdot v_e>_ed\mu_e=4\pi
\end{equation}
because the origin is on the sphere $S^2(a)$. Since

\begin{eqnarray}
  \lim_{n\rightarrow\infty}\int_{\Sigma_n\cap
  B^c_{sH^{-1}}(0)}\frac{2}{r^3}<x\cdot v_e>_ed\mu_e=\int_{S^2(a)\cap B^c_s(0)}\frac{2}{r^3}<x\cdot v_e>_ed\mu_e
\end{eqnarray}
\begin{eqnarray}
  \lim_{n\rightarrow\infty}\int_{\Sigma_n\cap B_{Kr_0}(0)}\frac{2}{r^3}<x\cdot
  v_e>_ed\mu_e=\int_{P\cap B_K(0)}\frac{2}{r^3}<x\cdot
  v_e>_ed\mu_e
\end{eqnarray}and

\begin{eqnarray}
  \int_P\frac{2}{r^3}<x\cdot v_e>_ed\mu_e=4\pi
\end{eqnarray}
then we have

\begin{eqnarray}\label{31}
  \lim_{s\rightarrow0,K\rightarrow\infty}\limsup_{n\rightarrow\infty}|\int_{\Sigma_n\cap (B_{sH^{-1}}\setminus B_{Kr_0})}\frac{2}{r^3}<x\cdot
  v_e>_ed\mu_e|=0
\end{eqnarray}

Now we want to prove that
\begin{eqnarray}\label{32}
  \lim_{s\rightarrow0,K\rightarrow\infty}\limsup_{n\rightarrow\infty}|\int_{\Sigma_n\cap (B_{sH^{-1}}\setminus B_{Kr_0})}\frac{2}{r^3}<x\cdot
  v_e>_e<v_e\cdot b>_ed\mu_e|=0
\end{eqnarray}

We use Lemma \ref{28} to get (\ref{32}) from (\ref{31}), but there
is a bit difference from \cite{QT}.

\begin{eqnarray}
 && \int_{\Sigma_n\cap (B_{sH^{-1}}\setminus
  B_{Kr_0})}\frac{2}{r^3}<x\cdot v_e>_e<v_e\cdot
  b>_ed\mu_e\nonumber\\
 && =<v_n\cdot b>_e\int_{\Sigma_n\cap (B_{sH^{-1}}\setminus
  B_{Kr_0})}\frac{2}{r^3}<v_e\cdot b>_ed\mu_e\nonumber\\
  &&+\int_{\Sigma_n\cap (B_{sH^{-1}}\setminus
  B_{Kr_0})}\frac{2}{r^3}<x\cdot v_e>_e<(v_e-v_n)\cdot
  b>_ed\mu_e
\end{eqnarray}
The first term will converge to 0. For the second term, we deal with
it in the cylinder coordinate in Section 4:

\begin{eqnarray}
  &&|\int_{\Sigma_n\cap (B_{sH^{-1}}\setminus
  B_{Kr_0})}\frac{2}{r^3}<x\cdot v_e>_e<(v_e-v_n)\cdot
  b>_ed\mu_e|\nonumber\\
  &&=|\sum_{j=1}^{l_n}\int_{A_{Kr_0e^{(j-1)L},Kr_0e^{jL}}}\frac{2}{r^3}<x\cdot v_e>_e<(v_e-v_n)\cdot
  b>_ed\mu_e|\nonumber\\
 &&\leq C\sum_{j=1}^{l_n}L\max_{I_j}|v_e-v_n|\nonumber\\
  &&=C\sum_{j=1}^{l_n/2}L\max_{I_j}|v_e-v_n|+C\sum_{j=l_n/2+1}^{l_n}L\max_{I_j}|v_e-v_n|
\end{eqnarray}

From Lemma \ref{28}
 \begin{eqnarray}
   &&CL\sum_{i=1}^{l_n/2}\sup_{I_i}|v-v_n|+CL\sum_{i=\frac{l_n}{2}+1}^{l_n}\sup_{I_i}|v-v_n|\nonumber\\
   &&\leq C(l_ne^{-\frac{1}{4}l_nL}+C)s+l_n^2r_0^{-\frac{1}{2}}
 \end{eqnarray}

But from the condition

\begin{equation}
  \lim_{n\rightarrow \infty}\frac{\log(r_1(\Sigma_n))}{r_{0}(\Sigma_n)^{1/4}}=0
\end{equation}
we know

\begin{equation}
  \lim_{n\rightarrow\infty}l_n^2r_0^{-\frac{1}{2}}=\lim_{n\rightarrow\infty}(\frac{L^{-1}(\log sH^{-1}-\log Kr_0)}{r_0^{\frac{1}{4}}})^2=0
\end{equation}
so (\ref{32}) holds.

Then
\begin{equation}
  0\leq-\frac{8}{3}\pi m(g)-\frac{16}{3}\pi m(g)+(1-\frac{2}{C})8\pi
  m(g)=-\frac{16}{C}\pi m(g)
\end{equation}
but $m(g)>0$ , this is a contradiction. So for the stable constant
mean curvature foliation there exists some constant $C>0$ such that
for any sphere $\Sigma$ in the foliation,

\begin{equation}
  \frac{r_0(\Sigma)}{r_1(\Sigma)}\geq C.
\end{equation}
Then the uniqueness follows from Theorem \ref{33}.

Proof of the Corollary \ref{34}. Suppose there is not such
$K(C,\beta)$, then we can find a sequence of  constant mean
curvature spheres $\Sigma_n$, with

\begin{eqnarray}
  \lim_{n\rightarrow\infty}r_0(\Sigma_n)=\infty&&\lim_{n\rightarrow\infty}\frac{\log(r_1)}{r_0^{\frac{1}{4}}}=0
\end{eqnarray}
and $\Sigma_n$ do not belong to the foliation. But from the argument
above we know this sequence satisfies
\begin{equation}
  \frac{r_0(\Sigma_n)}{r_1(\Sigma_n)}\geq C.
\end{equation}
So when $n$ is sufficiently large, $\Sigma_n$ must belong to the
foliation, which ends the proof.

\end{document}